\documentclass[12pt]{amsart}
\usepackage{color}
\usepackage{epsfig}
\usepackage{verbatim}

\usepackage{algorithm}
\usepackage{algpseudocode}
\usepackage{amssymb}

\newtheorem{theorem}{Theorem}[section]

\newtheorem{lemma}[theorem]{Lemma}

\begin{document}

\title {Constructing Minimally 3-connected graphs}
\maketitle 
 
\begin{center}
Jo\~{a}o Paulo Costalonga \footnote{Departamento de Matem\'atica 
Universidade Federal Do Esp\'irito,
Av. Fernando Ferrari 514,
 Vit\'oria, Esp\'irito Santo, 29075-910, Brazil,
joaocostalonga@gmail.com.}, 
R. J. Kingan \footnote{Bloomberg LP., 
120 Park Avenue,
New York, NY 10165,
rkingan@bloomberg.net.}, 
S. R. Kingan \footnote{Department of Mathematics,
Brooklyn College, City University of New York,
 Brooklyn, NY 11210,
skingan@brooklyn.cuny.edu.}  
\end{center}
   
\begin{abstract}
A $3$-connected graph is minimally 3-connected if removal of any edge destroys 3-connectivity. We present an algorithm for constructing minimally 3-connected graphs based on the results in (Dawes, JCTB 40, 159-168, 1986) using two operations: adding an edge between non-adjacent vertices and splitting a vertex. To test sets of vertices and edges for 3-compatibility, which~depends on the cycles of the graph, we develop a method for obtaining the cycles of $G'$ from the cycles of $G$, where $G'$ is obtained from $G$ by one of the two operations above.  We eliminate isomorphic duplicates using certificates generated by McKay's isomorphism checker nauty. The algorithm consecutively constructs the non-isomorphic minimally 3-connected graphs with $n$ vertices and $m$ edges from the non-isomorphic minimally 3-connected graphs with $n-1$ vertices and $m-2$ edges,  $n-1$ vertices and $m-3$ edges, and $n-2$ vertices and $m-3$ edges.
\end{abstract}

\section{Introduction}\label{sec1}
 \bigskip

\begin{figure}[H]
\centering
\includegraphics[width=3.8in]{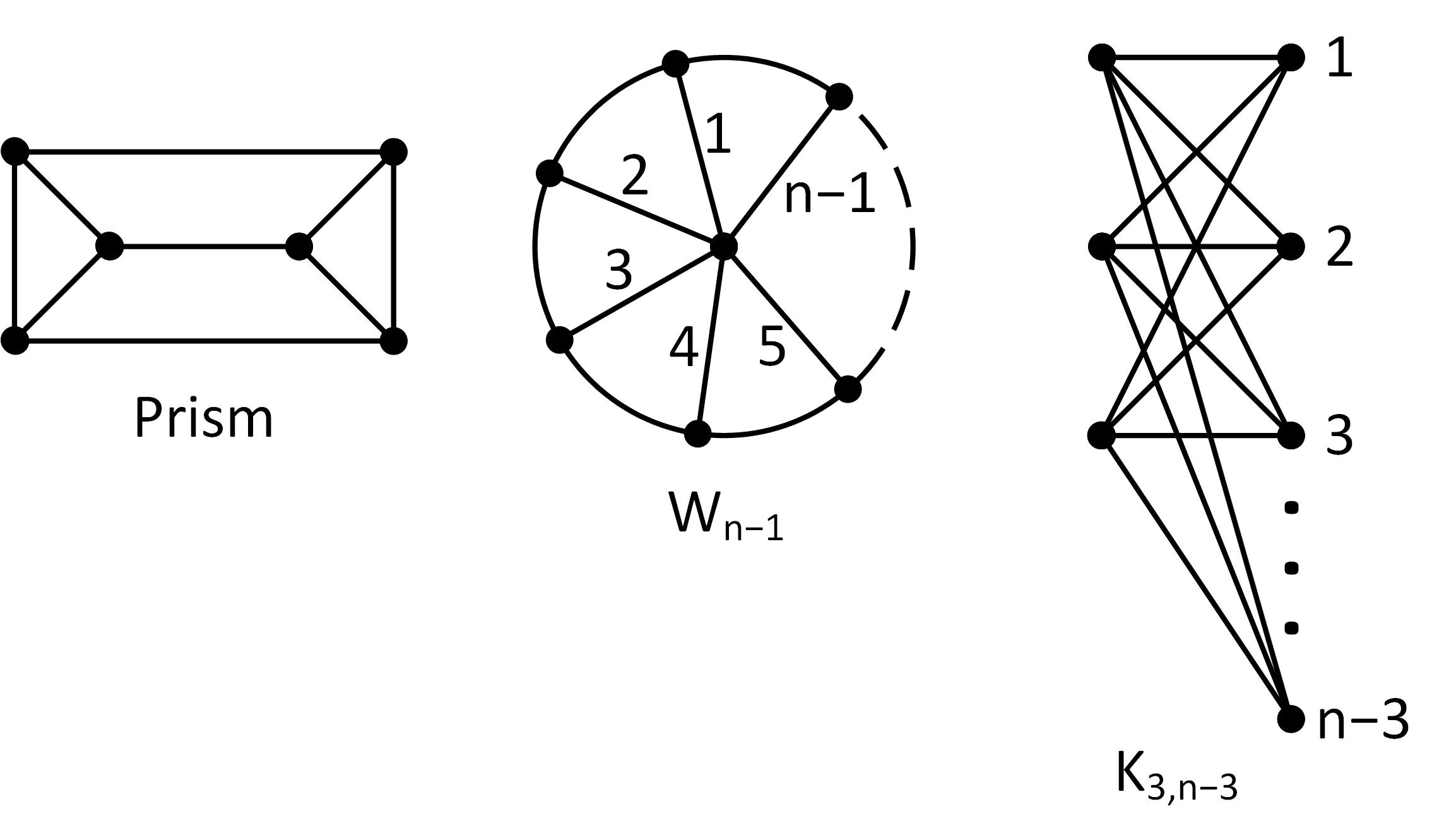}
\caption{ The prism graph, $W_{n-1}$ for $n \ge 4$, and $K_{3, n-3}$ for $n \ge 6$. \label{Exprism}} 
\end{figure}

In this paper, we present an algorithm for consecutively generating minimally 3-connected graphs, beginning with the prism graph, with the exception of two families. The two exceptional families are the wheel graph with $n$ vertices and $n-1$ spokes denoted by $W_{n-1}$ for $n \ge 4$ and the complete bipartite graph with 3 vertices in one class and $n-3$ vertices in the other class denoted by $K_{3, n-3}$ for $n \ge 6$. See Figure \ref{Exprism}.

\bigskip

\begin{figure}[H]
\centering
\includegraphics[width=2.9in]{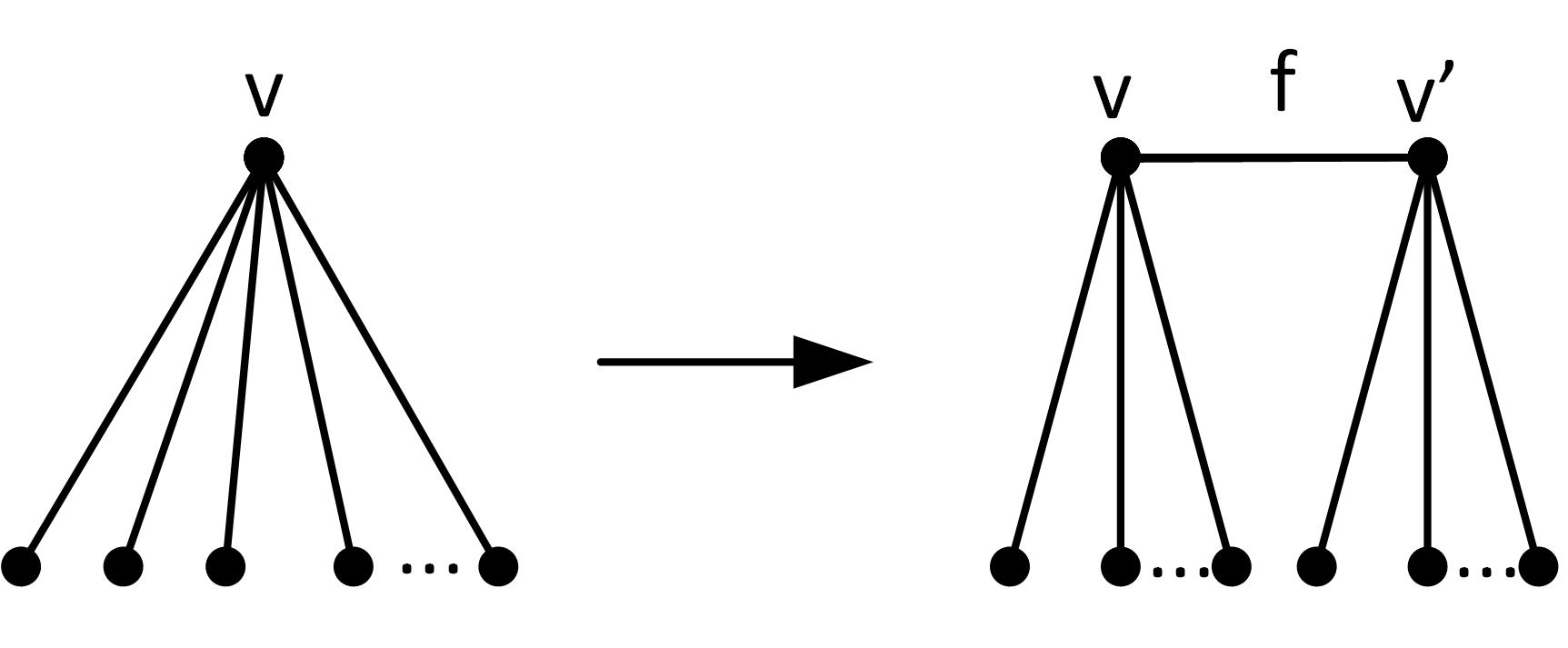}
\caption{ Tutte's vertex split operation. \label{VertexSplit}}
\end{figure}

Tutte proved that a simple graph is $3$-connected if and only if it is a wheel or is obtained from a wheel by adding edges between non-adjacent vertices and splitting vertices [\ref{Tutte1961}]. The vertex split operation is illustrated in Figure \ref{VertexSplit}. At each stage the graph obtained is 3-connected.

\begin{figure}[H]
\centering
\includegraphics[width=2.45in]{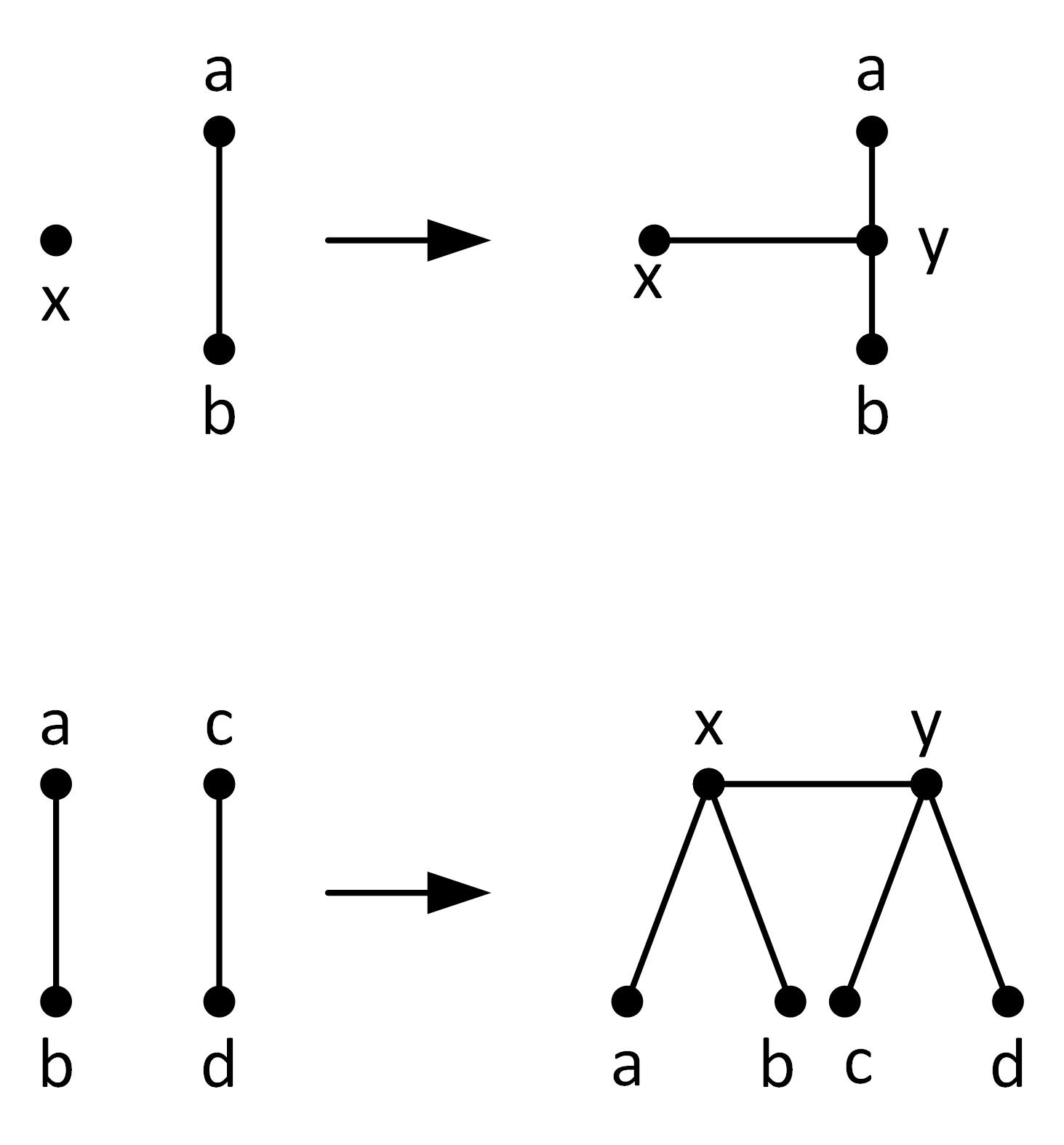}
\caption{ Bridging two edges and bridging a vertex and an edge. \label{BG-1and2}}
\end{figure}

A cubic graph is a graph whose vertices have degree 3.   Let $G$ and $H$ be 3-connected cubic graphs such that $G\not\cong W_3$ and $H$ is a minor of $G$.  A pair of distinct edges is \textit{bridged} if they are subdivided by vertices $x$ and $y$, respectively, forming paths of length 2, and $x$ and $y$ are joined by an edge. This~operation is explained in detail in Section \ref{Section2} and illustrated in Figure \ref{BG-1and2}.  Tutte also proved that $G$ can be obtained from $H$ by repeatedly bridging edges. At each stage the graph obtained remains 3-connected and cubic [\ref{Tutte1967}]. Observe that this operation is equivalent to adding an edge $e=uv$ to a cubic graph and splitting $u$ and splitting $v$. This gives an easy way of consecutively constructing all 3-connected cubic graphs on $n$ vertices for even $n$. Surprisingly the entry for the number of 3-connected cubic graphs in the Online Encyclopedia of Integer Sequences (sequence A204198) has entries only up to $n=14$. We were able to quickly obtain such graphs up to $n=20$.

Tutte's result and our algorithm based on it suggested that a similar result and algorithm may be obtainable for the much larger class of minimally 3-connected graphs.

\medskip
Dawes  gave a necessary and sufficient characterization for the construction of minimally 3-connected graphs starting with $W_3$. To do this he needed three operations one of which is the above operation where two distinct edges are bridged. Let $e=ab$ be an edge and $x\neq a, b$ be a vertex. A vertex and an edge are \textit{bridged} if a new vertex is placed on edge $e$ and linked to $x$. Dawes proved that starting with $W_3$ the class of minimally 3-connected graphs can be constructed by bridging a vertex and an edge, bridging two edges, or by adding a degree 3 vertex in the manner Dawes specified using what he called ``3-compatible sets'' as explained in Section \ref{Section2}.  Following the above approach for cubic graphs we were able to translate Dawes' operations to edge additions and vertex splits and develop an algorithm that consecutively constructs minimally 3-connected graphs from smaller minimally 3-connected graphs.

\medskip
First, we prove exactly how Dawes' operations can be translated to edge additions and vertex~splits. Second, we prove a cycle propagation result. Although obtaining the set of cycles of a graph is NP-complete in general, we can take advantage of the fact that we are beginning with a~fixed cubic initial graph, the prism graph. We develop methods for constructing the set of cycles for a graph $G'$ obtained from a graph $G$ by edge additions and vertex splits, and Dawes specifications on 3-compatible sets.    Let $n$ be the number of vertices in $G$ and let $c$ be the number of cycles of $G$. We prove that the set of cycles of $G'$ can be obtained from the set of cycles of $G$ by a method with complexity~$\mathcal O(c^2 n)$.   Third, we prove that if $G$ is a minimally 3-connected graph that is not $W_{n-1}$ for $n\ge 4$ or $K_{3, n-3}$ for $n\ge 6$, then  $G$ must have a prism minor, for $n \ge 7$, and $G$ can be obtained from a smaller minimally $3$-connected graph $G'$ such that $|E(G)| - |E(G')| \le 3$ using edge additions and vertex splits and Dawes specifications on 3-compatible sets.

\medskip
 
We present an algorithm based on the above results that consecutively constructs the non-isomorphic minimally 3-connected graphs with $n$ vertices and $m$ edges from the non-isomorphic minimally 3-connected graphs with $n-1$ vertices and $m-2$ edges,  $n-1$ vertices and $m-3$ edges, and~$n-2$ vertices and $m-3$ edges. This formulation also allows us to determine worst-case complexity for processing a single graph; namely $\mathcal O(c^2n^3)$, which includes the complexity of cycle propagation mentioned above.

\medskip

There has been a significant amount of work done on identifying efficient algorithms for certifying $3$-connectivity of graphs. Hopcroft and Tarjan published a linear-time algorithm for testing $3$-connectivity [\ref{Hopcroft1973}]. Schmidt  extended this result by identifying a certifying algorithm for checking $3$-connectivity in linear time  [\ref{Schmidt2011}]. The perspective of this paper is somewhat different. Instead of checking an existing graph to determine whether it is minimally $3$-connected, we seek to construct graphs from the prism using a procedure that generates only minimally $3$-connected graphs.
The~algorithm presented in this paper is the first to generate exclusively minimally 3-connected graphs from smaller minimally 3-connected graphs.

\bigskip

%%%%%%%%%%%%%%%%%%%%%%%%%%%%%%%

\section {\bf Terminology, Previous Results, and Outline of the Paper} \label{Section2}
\bigskip

We begin with the terminology used in the rest of the paper. Since graphs used in the paper are not necessarily simple, when they are it will be specified. Let $G$ be a graph and $e=uv$ be an edge with end vertices $u$ and $v$. The graph with edge $e$ deleted is called an {\it edge-deletion} and is denoted by $G\backslash e$ or~$G\backslash uv$. When deleting edge $e$, the end vertices $u$ and $v$ remain. 
To contract edge $e$, collapse the edge by identifing the end vertices $u$ and $v$ as one vertex, and delete the resulting loop. The graph with edge $e$~contracted is called an {\it edge-contraction} and denoted by  $G / e$. 
A graph $H$ is   a {\it minor} of a~graph $G$ if $H$ can be obtained from $G$ by deleting edges (and any isolated vertices formed as a result) and contracting~edges. 
We write $H=G\backslash X/Y$, where $X$ is the set of edges deleted and $Y$ is the set of edges~contracted. 

\medskip

A {\it triangle} is a set of three edges in a cycle and a {\it triad} is a set of three edges incident to a degree 3 vertex.   A graph is {\it $3$-connected} if at least 3 vertices must be removed to disconnect the graph. In~a~3-connected graph $G$, an edge $e$ is {\it deletable} if $G\backslash e$ remains 3-connected. A 3-connected graph with no deletable edges is called {\it minimally $3$-connected}.

\medskip

There are multiple ways that deleting an edge in a minimally 3-connected graph $G$ can destroy connectivity. One obvious way is when $G$ has a degree 3 vertex $v$ and deleting one of the edges incident to $v$ results in a 2-connected graph that is not 3-connected. 
Halin proved  that a minimally 3-connected graph has at least one triad  [\ref{Halin1969}]. We exploit this property to develop a construction theorem for minimally 3-connected graphs.

\medskip

The operation that reverses edge-deletion is edge addition. 
A simple graph $G$ with an edge $e=uv$ added between non-adjacent vertices is called an {\it edge addition} of $G$ and   denoted by $G+e$ or $G+uv$. 

\bigskip
The operation that reverses edge-contraction is called a vertex split of $G$. To {\it split} a vertex $v$ with $deg_G(v) \ge 4$, first divide $N_G(v)$ into two disjoint sets $S$ and $T$, both of size at least 2. Then replace $v$ with two distinct vertices $v$ and $v'$, join them by a new edge $f=vv'$, and join each neighbor of $v$ in $S$ to $v$ and each neighbor in $T$ to $v'$. The resulting graph is called a {\it vertex split}  of $G$ and is denoted by $K = G\circ_{S, T} f$. In other words $N_G(v)$ is partitioned into two sets $S$ and $T$, and in $K$, $N_{K}(v) = S \cup v'$ and $N_{K}(v') = T \cup v$.  Observe that $deg_K(v) \ge 3$ and $deg_K(v') \ge 3$.

\bigskip

We can get a different graph depending on the assignment of neighbors of $v$ in $G$ to $v$ and $v'$ in the vertex split; hence the sets $S$ and $T$ in the notation. With a slight abuse of notation, we can say $G \circ f$, as each vertex split is described with a particular assignment of neighbors of $v$ in $G$ to $v$ and $v'$. When performing a vertex split, we will think of $v'$ as the new vertex that gets added and $f=vv'$ as the new edge that gets added. Observe that  if $G$ is 3-connected, then edge additions and vertex splits remain 3-connected. The degree condition $deg_G(v) \ge 4$ is not necessary for an arbitrary vertex split, but required to preserve 3-connectivity. Figure \ref{VertexSplit} shows the vertex split operation. 

\medskip

In 1961 Tutte proved that
a simple graph  is $3$-connected if and only if it is a wheel or is obtained from a wheel by a finite sequence of edge additions or vertex splits. This result is known as Tutte's Wheels Theorem  [\ref{Tutte1961}].  

\medskip

\begin{theorem} \label{TutteWheelsTheorem} {\bf (Tutte, 1961)} Let $G$ be a simple graph that is not a wheel. Then $G$ is $3$-connected if and only if $G$ can be constructed from a wheel minor by a finite sequence of edge additions or vertex splits. \qed
\end{theorem}  

\medskip
The next result we need  is Dirac's characterization of 3-connected graphs without a prism minor  [\ref{Dirac1963}].  The graphs $K'_{3,n-3} $, $K''_{3,n-3} $, and $K'''_{3,n-3}$ are obtained from the complete bipartite graph $K_{3, n-3}$ (shown in Figure \ref{Exprism}) with one, two, or three edges, respectively, joining the three vertices in one class.

\medskip
\begin{theorem} \label{DiracTheorem} {\bf (Dirac, 1963)}  A simple $3$-connected graph $G$ has no prism-minor if and only if $G$ is isomorphic to $K_5\backslash e$, $K_5$, $W_{n-1}$, for $n\ge 4$, $K_{3,n-3}$,  $K'_{3,n-3} $, $K''_{3,n-3}$, or $K'''_{3,n-3}$, for $n\ge 6$.  \qed
\end{theorem}  

\medskip

Theorem \ref{DiracTheorem} implies that there are only two infinite families of minimally 3-connected graphs without a prism-minor, namely $W_{n-1}$ for $n\ge 4$ and $K_{3,n-3}$ for $n\ge 6$. Thus, we may focus on constructing minimally 3-connected graphs with a prism minor.  

\medskip

Next, Halin proved that minimally 3-connected graphs are sparse in the sense that there is a linear bound on the number of edges in terms of the number of vertices   [\ref{Halin1969}].

\medskip
\begin{theorem} \label{HalinTheorem} {\bf (Halin, 1969)}  Let $G$ be a minimally $3$-connected graph on $n\ge 8$  vertices. Then   $|E(G)|\le 3n-9$. Moreover, $|E(G)|= 3n-9$ if and only if $G\cong K_{3, n-3}$. \qed
\end{theorem}  

\medskip

In 1969 Barnette and Gr\"{u}nbaum defined two operations based on subdivisions and gave an alternative construction theorem for 3-connected graphs   [\ref{BarnetteGrunbaum1969}]. A {\it subdivision} of $G$ is obtained from $G$ by replacing an edge by a path of length at least 2.
Let $G$ be a 3-connected graph. The first Barnette and Gr\"{u}nbaum operation is defined as follows: Subdivide an edge $ab$  by vertex $y$ and add edge $xy$ for a vertex $x\neq a, b$. This is what we called ``bridging a vertex and an edge'' in Section \ref{sec1}. The~second Barnette and Gr\"{u}nbaum operation is defined as follows: Subdivide two distinct edges $ab$ and $cd$, by~vertices $x$ and $y$, respectively, and add edge $xy$. This is what we called ``bridging two edges'' in~Section \ref{sec1}. Observe that these operations, illustrated in Figure \ref{BG-1and2}, preserve 3-connectivity.     

\medskip

\begin{theorem} \label{BarnetteGrunbaumTheorem} {\bf (Barnette and Gr\"{u}nbaum, 1968)} Let $G$ be a simple graph such that $G\not\cong W_3$. Then $G$ is $3$-connected if and only if $G$ can be constructed from $W_3$ by a finite sequence of edge additions, bridging a vertex and an edge, or bridging two edges. \qed
\end{theorem}

\medskip

\begin{figure}[H]
\centering
\includegraphics[width=2.25in]{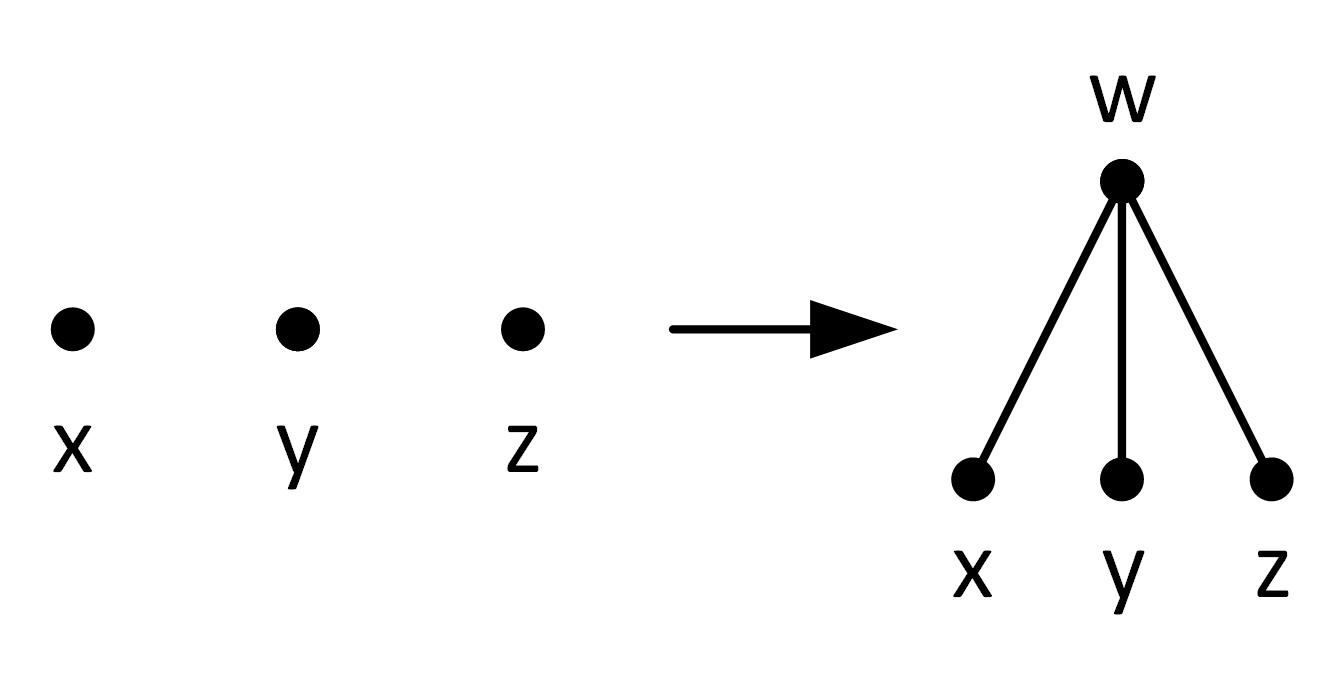}
\caption{ Dawes' third operation. \label{D-3}} 
\end{figure}

In 1986, Dawes gave a necessary and sufficient characterization for the construction of minimally 3-connected graphs starting with $W_3$. He used the two Barnett and Gr\"{u}nbaum operations (bridging an edge and bridging a vertex and an edge) and a new operation, shown in Figure \ref{D-3}, that he defined as follows: select three  distinct vertices $x, y, z$  in the graph and link all three to a new vertex $w$ by adding three new edges $xw$, $yw$, and $zw$. Observe that this new operation also preserves 3-connectivity. We will call this operation ``adding a degree 3 vertex'' or in matroid language ``adding a triad'' since a triad is a set of three edges incident to a degree 3 vertex. 
Using these three operations, Dawes gave a necessary and sufficient condition for the construction of minimally 3-connected graphs. 

\begin{figure}[H]
\centering
\includegraphics[width=1.45in]{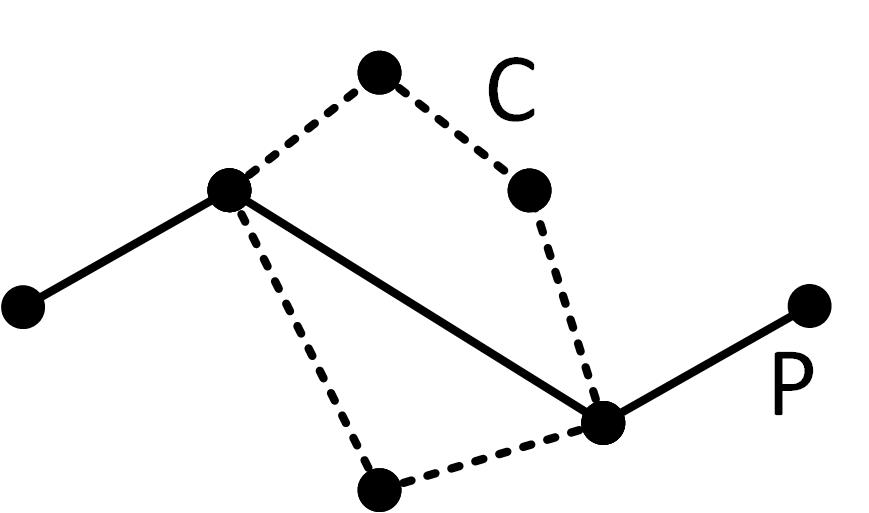}
\caption{ A chording path. \label{chord}}
\end{figure}

 Dawes thought of the three operations, bridging edges, bridging a vertex and an edge, and the third operation as acting on, respectively,
a vertex and an edge, two edges, and three vertices. 
Dawes showed that if one begins with a minimally 3-connected graph and applies one of these operations, the~resulting graph will also be minimally 3-connected if and only if certain conditions are met.
Let~$C$ be a~cycle in a graph $G$. A {\it chord} of $C$ is an edge $e \not\in C$ that links two vertices in $C$. A {\it chording path} $P$ for a~cycle $C$ is a path that has a chord $e$ of $C$ in it and intersects $C$ only in the end vertices of $e$. In~particular, none of the edges of $C$ can be in the path. See Figure \ref{chord}.

When applying the three operations listed above, Dawes defined conditions on the set of vertices and/or edges being acted upon that guarantee that the resulting graph will be minimally $3$-connected. A set $S$ of vertices and/or edges in a graph $G$ is {\it $3$-compatible} if it conforms to one of the following three~types:
\begin{enumerate}
\item $S = \{x, ab\}$, where $x$ is a vertex of $G$, $ab$ is an edge of $G$, $x \neq a, b$ and no $xa$-path or $xb$-path is a~chording path of $G - ab$;
\item $S = \{ab, cd\}$, where $ab$ and $cd$ are distinct edges of $G$, though possibly adjacent, and no $ac$-, $bc$-, $ad$- or $bd$-path is a chording path of $G - \{ab, cd\}$; or
\item $S = \{x, y, z\}$, where $x$, $y$, and $z$ are distinct vertices of $G$ and no $xy$-, $xz$- or $yz$-path is a chording path of $G$. Please note that if $G$ is $3$-connected, then $x$, $y$, and $z$ must be pairwise non-adjacent if $\{x, y, z\}$ is 3-compatible.
\end{enumerate}
For convenience in the descriptions to follow, we will use D1, D2, and D3 to refer to bridging a~vertex and an edge, bridging two edges, and adding a degree 3 vertex, respectively.
Dawes proved that if one of the operations D1, D2, or D3 is applied to a minimally $3$-connected graph, then the result is minimally $3$-connected if and only if the operation is applied to a 3-compatible set  [\ref{Dawes1986}].

\medskip
\begin{theorem} \label{DawesTheorem} {\bf (Dawes, 1986a)} Let $H$ be a minimally $3$-connected graph. Let $G$ be constructed from $H$ by applying D1, D2, or D3 to a set $S$ of edges and/or vertices of $H$. Then $G$ is minimally $3$-connected if and only if $S$ is a $3$-compatible set in $H$. \qed
\end{theorem}

\medskip

Dawes also proved that, with the exception of $W_3$, every minimally 3-connected graph  can be obtained by applying D1, D2, or D3 to a 3-compatible set in a smaller minimally 3-connected graph.

\medskip

\begin{theorem} \label{DawesTheorem2} {\bf (Dawes, 1986b)} Let $G$ be a simple graph such that $G\not \cong W_3$. Then $G$ is minimally $3$-connected if and only if there exists a minimally $3$-connected graph $G'$, $|E(G')| < |E(G)|$ such that $G$ can be constructed by applying one of D1, D2, or D3 to a $3$-compatible set in $G'$. \qed
\end{theorem}

\medskip

 The next result is the Strong Splitter Theorem  [\ref{KinganLemos2014}]. The rank of a graph, denoted by $r(G)$, is the size of a spanning tree. If $G$ has $n$ vertices, then $r(G)=n-1$.

\medskip
\begin{theorem}\label{strong-splitter-theorem}
Suppose $G$ and $H$ are simple $3$-connected graphs such that $G$ has a proper $H$-minor, $G$ is not a~wheel, and $H\not\cong W_3$. Let $j=r(G)-r(H)$. Then there is a sequence of $3$-connected graphs $G_0,G_1,\dots,G_t$~such that $G_0\cong H$, $G_t=G$, and $G_{i-1}$ is a minor of $G_i$ such that:
\begin{enumerate}
\item[(i)] For $1\le i \le j$, $r(G_i)-r(G_{i-1})=1$ and $|E(G_i)|-|E(G_{i-1})|\le 3$; and
\item[(ii)] For $j< i\le t$, $r(G_i)=r(G)$ and $|E(G_i)|-|E(G_{i-1})|=1$.
\end{enumerate}
Moreover, when $|E(G_i)|-|E(G_{i-1})|=3$, for $1\le i\le j$,  $E(G_i)-E(G_{i-1})$ is a triad of $G_i$. \qed
\end{theorem}

\medskip

Our goal is to generate all minimally 3-connected graphs with $n$ vertices and $m$ edges, for various values of $n$ and $m$ by repeatedly applying operations D1, D2, and D3 to input graphs after checking the input sets for 3-compatibility. The process needs to be \textit{correct}, in that it only generates minimally 3-connected graphs, \textit{exhaustive}, in that it generates all minimally 3-connected graphs, and \textit{isomorph-free}, in that no two graphs generated by the algorithm should be isomorphic to each other.

\medskip

By Theorem \ref{DawesTheorem}, in order for our method to be correct it needs to verify that a set of edges and/or vertices is $3$-compatible before applying operation D1, D2, or D3.  In Section \ref{correctness-section}, we present two of the three new theorems in this paper. The first new result expresses operations D1, D2, and D3 in terms of edge additions and vertex splits. The second new result gives an algorithm  for the efficient propagation of the list of cycles of a graph from a smaller graph when performing edge additions and vertex splits. We call it the ``Cycle Propagation Algorithm.'' Together, these two results establish correctness of the method. In Section \ref{cycle-propagation-section} we provide details of the implementation of the Cycle Propagation Algorithm.

\medskip

In Section \ref{InfiniteBookshelf-section} we present the algorithm for generating minimally 3-connected graphs using an~``infinite bookshelf'' approach to the removal of isomorphic duplicates by lists. Specifically, we show how we can efficiently remove isomorphic graphs from the list of generated graphs by restructuring the operations into atomic steps and computing only graphs with fixed edge and vertex counts in~batches.

\medskip

In Section \ref{exhaustive-section} we show that the ``Infinite Bookshelf Algorithm'' described in Section \ref{InfiniteBookshelf-section}  is exhaustive by showing that all minimally 3-connected graphs with the exception of two infinite families, $W_{n-1}$ and $K_{3, n-3}$, can be obtained from the prism graph by applying operations D1, D2, and D3.  This is the third new theorem in the paper.
\bigskip

%%%%%%%%%%%%%%%%%%%%%%%%%%%%%%%%%%%%%%%%%%%%%%

\section{\bf Results Establishing Correctness of the Algorithm}\label{correctness-section} 

\bigskip

In this section, we present two results that establish that our algorithm is correct; that is, that it produces only minimally 3-connected graphs. 

\medskip

According to Theorem \ref{DawesTheorem}, when operation D1, D2, or D3 is applied to a set $S$ of edges and/or vertices in a minimally 3-connected graph, the result is minimally 3-connected if and only if $S$ is 3-compatible. To check whether a set is 3-compatible, we need to be able to check whether chording paths exist between pairs of vertices. To check for chording paths, we need to know the cycles of the graph. Since enumerating the cycles of a graph is an NP-complete problem, we would like to avoid it by determining the list of cycles of a graph generated using D1, D2, or D3 from the cycles of the graph it was generated from. 

\medskip
To determine the cycles of a graph produced by D1, D2, or D3, we need to break the operations down into smaller ``atomic'' operations. The first theorem in this section, Theorem \ref{maintheorem-1}, expresses operations D1, D2, and D3 in terms of edge additions and vertex splits. The second theorem in this section, Theorem \ref{cycle-prop-theorem}, provides bounds on the complexity of a procedure to identify the cycles of a~graph generated through operations D1, D2, and D3 from the cycles of the original graph. The~second theorem relies on two key lemmas which show how cycles can be propagated through edge additions and vertex splits. We refer to these lemmas multiple times in the rest of the paper.

\medskip

\begin{theorem}\label{maintheorem-1} Let $G$ be a simple $3$-connected graph. Operations D1, D2, and D3 can be expressed as a sequence of edge additions and vertex splits. Specifically:

\begin{itemize}

\item[(a)] D1 applied to a vertex $x$ and an edge $ab$ in $G$ to create a new edge $xy$ can be expressed as $(G + e) \circ f$, where $e = xa$ and $f = ay$;

\item[(b)] D2 applied to two edges $ab$ and $cd$ in $G$ to create a new edge $xy$ can be expressed as $(G + e) \circ \{f_c, f_b\}$, where $e = bc$, $f_b = xb$ and $f_c = cy$; and

\item[(c)] D3 applied to vertices $x$, $y$ and $z$ in $G$ to create a new vertex $w$ and edges $xw$, $yw$ and $zw$ can be expressed as $(G + \{e_1, e_2\}) \circ f$, where $e_1 = yx$, $e_2 = zx$ and $f = xw$.
\end{itemize}
\end{theorem}

\begin{figure}[H]
\centering
\includegraphics[width=5.05in]{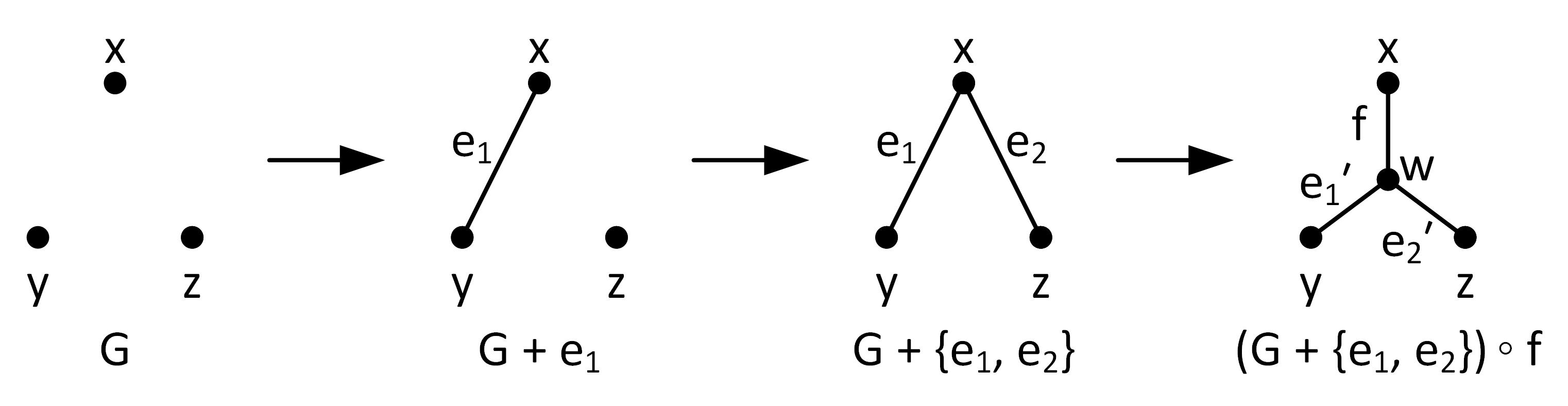}
\caption{Operation D3 interpreted as two edge additions followed by a vertex split.\label{D-3-equivalent}}
\end{figure}

\begin{proof}
 Operation D1 requires a vertex $x$ and a nonincident edge $ab$. The operation is performed by subdividing edge $ab$ by vertex $y$, and adding edge $xy$. We may also interpret this operation as adding an edge $e = xa$, and then splitting vertex $a$ in such a way that $y$ is the new vertex adjacent to $x$ and~$b$, and the new edge $f = ya$, as shown in Figure \ref{BG-1-equivalent}. In the process, edge $e=xa$ is replaced with a new edge $e'=xy$ and edge $ab$ is replaced with a new edge $yb$. Following this interpretation, the resulting graph is~$(G + e) \circ f$.

\medskip

\begin{figure}[H]
\centering
\includegraphics[width=3.45in]{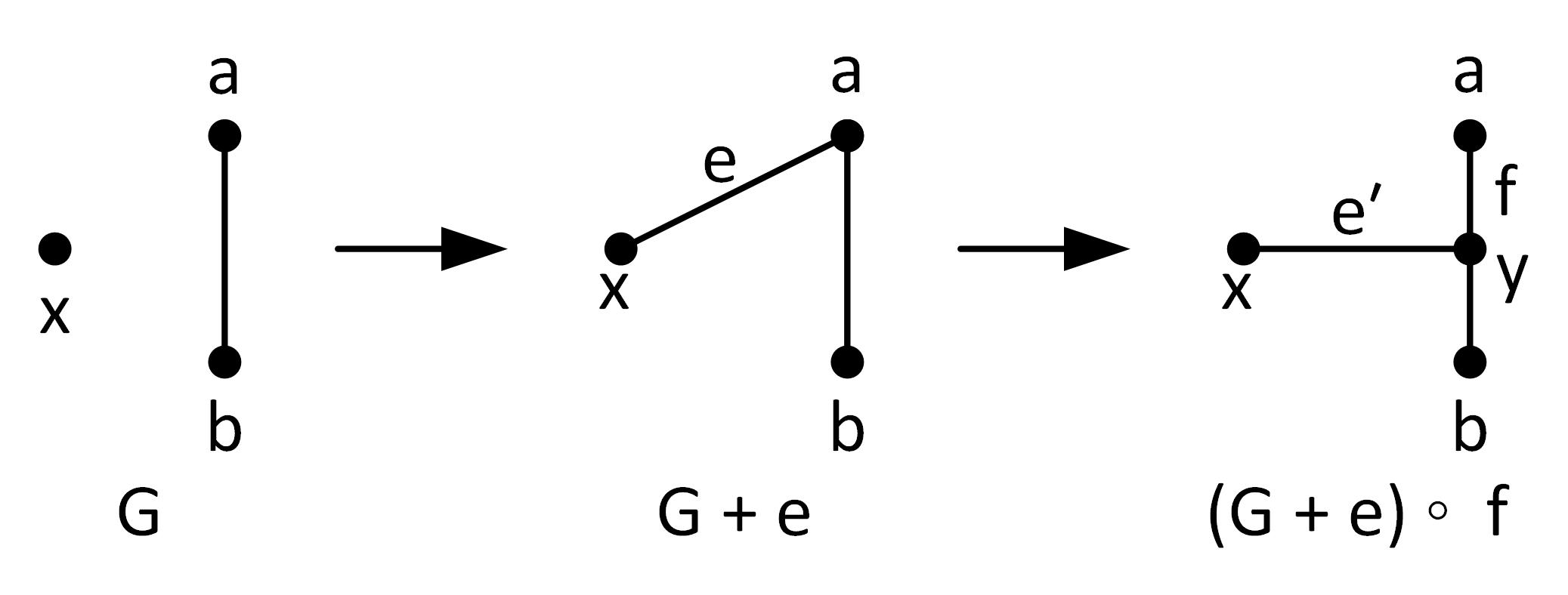}
\caption{Operation D1 interpreted as an edge addition followed by a vertex split.\label{BG-1-equivalent}}
\end{figure}

Operation D2 requires two distinct edges $ab$ and $cd$, and is performed by subdividing both edges and adding a new edge connecting the two vertices. We may interpret this operation using the following steps, illustrated in Figure \ref{BG-2-equivalent}:
\begin{itemize}
\item[(i)] Add an edge $e=bc$;
\item[(ii)] split the vertex $c$ in such a way that $y$ is the new vertex adjacent to $b$ and $d$, and the new edge $f_c = yc$; and
\item[(iii)] split the vertex $b$ in such a way that $x$ is the new vertex adjacent to $a$ and $y$, and the new edge $f_b = bx$.
\end{itemize}
\noindent In step (ii), edge $e=bc$ is replaced with a new edge $e'=yb$ and edge $cd$ is replaced with a new edge $yd$. In step (iii), edge $e'=yb$ is replaced with a new edge $e''=xy$ and $ab$ is replaced with a new edge $xa$. Following this interpretation, the resulting graph is $(G + e) \circ \{f_c, f_b\}$.
\medskip

\begin{figure}[H]
\centering
\includegraphics[width=5.05in]{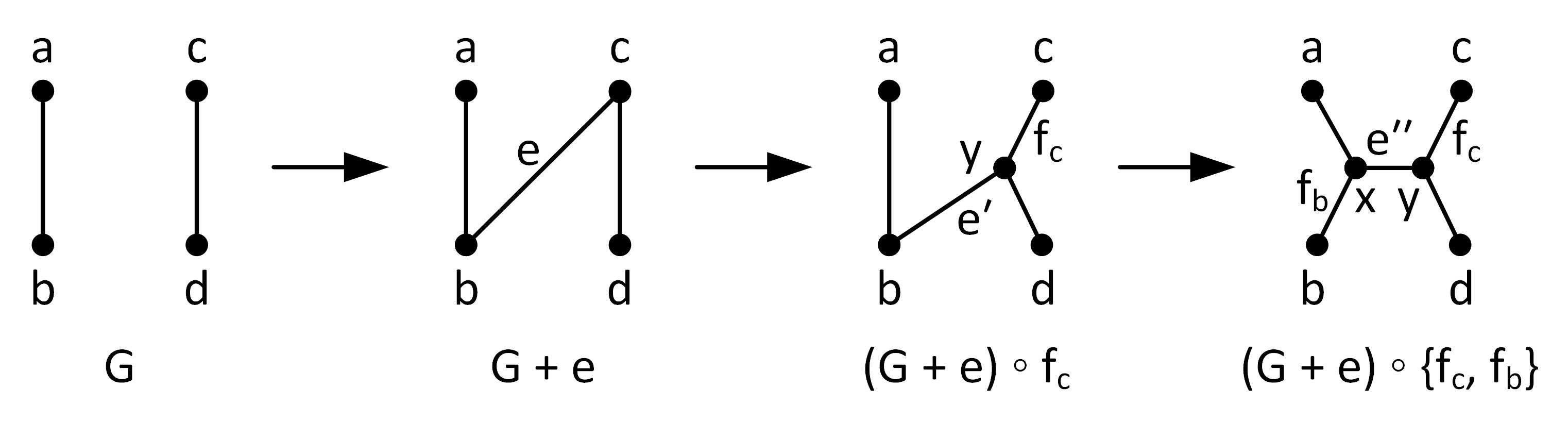}
\caption{Operation D2 interpreted as an edge addition followed by two vertex splits.\label{BG-2-equivalent}}
\end{figure}

Operation D3 requires three vertices $x$, $y$, and $z$. The operation is performed by adding a new vertex $w$ and edges $xw$, $yw$, and $zw$. We may interpret this operation as adding one edge $e_1 = xy$, adding a second edge $e_2 = xz$, and then splitting the vertex $x$ in such a way that $w$ is the new vertex adjacent to $y$ and $z$, and the new edge $f = xw$. This is illustrated in Figure \ref{D-3-equivalent}. Following this interpretation, the resulting graph is $(G + \{e_1, e_2\}) \circ f$. \end{proof}

\medskip

In Theorem \ref{maintheorem-1}, it is possible that the initially added edge in each of the sequences above is a~parallel~edge; however we will see in Section \ref{exhaustive-section} that we can avoid adding parallel edges by selecting our initial ``seed'' graph carefully. 

\medskip

Consider the function \textproc{HasChordingPath}$(G, a, b, K)$, where $G$ is a graph, $a$ and $b$ are vertices in $G$ and $K$ is a set of edges, whose value is True if there is a chording path from $a$ to $b$ in $G \backslash K$, and False otherwise. To efficiently determine whether $S$ is 3-compatible, whether $S$ is a set consisting of a vertex and an edge, two edges, or three vertices, we need to be able to evaluate \textproc{HasChordingPath}. To evaluate this function, we need to check all paths from $a$ to $b$ for chording edges, which in turn requires knowing the cycles of $G \backslash K$. The second theorem in this section establishes a bound on the complexity of obtaining cycles of a graph from cycles of a smaller graph. The proof consists of  two lemmas, interesting in their own right, and a short argument.

\medskip
Using Theorem \ref{maintheorem-1}, we can propagate the list of cycles of a graph through operations D1, D2, and~D3 if it is possible to determine the cycles of a graph $G'$ obtained from a graph $G$ by:
\begin{itemize}
\item Adding an edge $uv$ betweeen two non-adjacent vertices $u$ and $v$; and
\item Splitting a vertex $v$ in $G$ to form a new vertex $v'$ of degree 3 that is incident to the new edge $f=vv'$ and two other edges.  
\end{itemize}

\medskip

The first lemma shows how the set of cycles can be propagated when an edge $uv$ is added betweeen two non-adjacent vertices $u$ and $v$.

\medskip

\begin{figure}[H]
\centering
\includegraphics[width=2.5in]{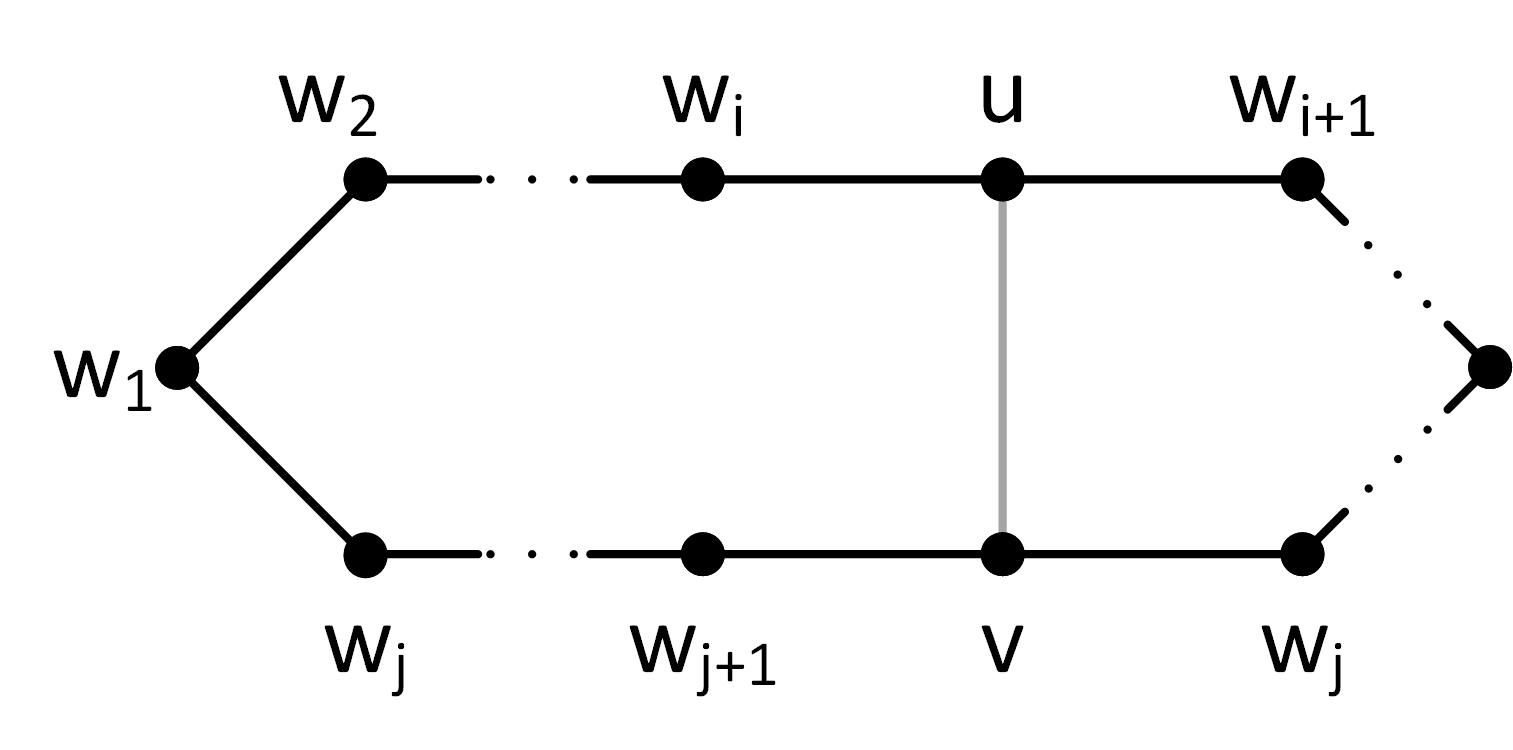}
\caption{Chording a cycle \label{cycle-chording-figure}}
\end{figure}

\begin{lemma} \label {CycleChordingLemma}  {\bf (Cycle Chording Lemma)} Let $G$ be a simple $2$-connected graph with $n$ vertices and let $\mathcal C(G)$ be the set of cycles of $G$. Let $G'$ be obtained from $G$ by adding an edge $uv$ between two non-adjacent vertices in $G$. Then the cycles of $G'$ consists of:
\begin{enumerate}
\item[(i)] $\mathcal C(G)$; and
\item[(ii)] Cycles $$w_1, \ldots, w_i, u, v, w_{j+1}, \ldots, w_k, w_1$$ and $$u, w_{i+1}, \ldots, w_j, v, u,$$ where $$w_1, \ldots, w_i, u, w_{i+1}, \ldots, w_j, v, w_{j+1}, \ldots, w_k, w_1$$ is a cycle in $G$ passing through $u$ and $v$, as shown in Figure \ref{cycle-chording-figure}.
\end{enumerate}

The complexity of determining the cycles of $G'$ from the cycles of $G$ is $\mathcal O(|\mathcal C(G)|n)$.
\end{lemma}
 \medskip

\begin{proof} We need only show that any cycle in $G'$ can be produced by (i) or (ii). Suppose $C$ is a cycle in~$G'$. If $C$ does not contain the edge $uv$ then $C$ must also be a cycle in $G$. Otherwise, the edges in $C$ other than $uv$ form a $u-v$ path $P_1$ in $G$. Since $G$ is $2$-connected, there is another edge-disjoint $u-v$ path $P_2$ in $G$. Paths $P_1$ and $P_2$ together form a cycle in $G$, and $C$ can be obtained from this cycle using the operation in (ii) above.  Finally, the complexity of determining the cycles of $G'$ from the cycles of $G$ is $\mathcal O(|\mathcal C(G)|n)$ because each cycle has to be traversed once and the maximum number of vertices in a~cycle is $n$. 
\end{proof}
\medskip
 
The graph $G$ in the statement of Lemma \bigskip {CycleChordingLemma} must be 2-connected.
It is easy to find a counterexample when $G$ is not $2$-connected; adding an edge to a graph containing a bridge may produce many cycles that are not obtainable from cycles in $G$ by Lemma \ref{CycleChordingLemma} (ii).

\medskip

Obtaining the cycles when a vertex $v$ is split to form a new vertex $v'$ of degree 3 that is incident to the new edge $f=vv'$ and two other edges is more complicated. For the purpose of identifying cycles, we regard a vertex split, where the new vertex has degree 3, as a sequence of two ``atomic'' operations. Let $v$ be a vertex in a graph $G$ of degree at least 4, and let $p$, $q$, $r$, and $s$ be four other vertices in $G$ adjacent to $v$. The following two steps describe a vertex split of $v$ in which $p$ and $q$ become adjacent to the new vertex and $r$ and $s$ remain adjacent to $v$:

\begin{enumerate}
\item Subdivide the edge joining $v$ and $p$, adding a new vertex $v'$.
\item Remove the edge $qv$ and replace it with a new edge $qv'$.
\end{enumerate}

\begin{figure}[H]
\centering
\includegraphics[width=4.15in]{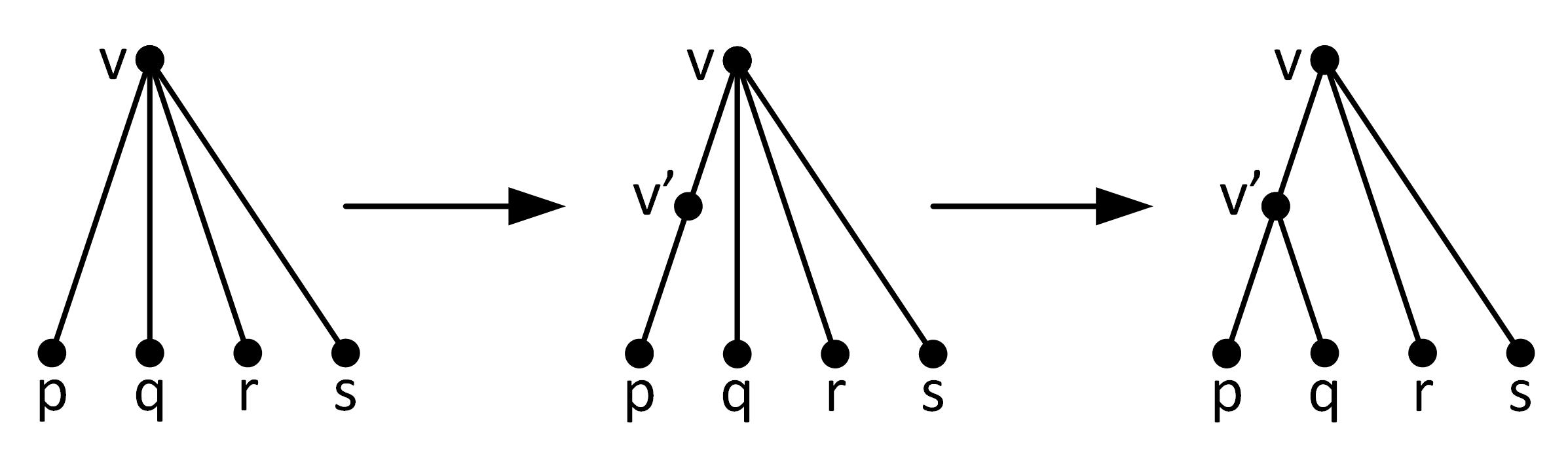}
\caption{ Rethinking a vertex split. \label{vertexsplitsforpropagatingcycles}}
\end{figure}

This is illustrated in Figure \ref{vertexsplitsforpropagatingcycles}. By thinking of the vertex split this way, if we start with the set of cycles of $G$, we can determine the set of cycles of $G \circ f$, where $G \circ f$ is obtained by splitting vertex $v$ to form a new vertex $v'$ of degree 3 that is incident to the new edge $f=vv'$ and two other edges. The~cycles of the graph resulting from step (1) above are simply the cycles of $G$, with any occurrence of the edge $pv$ replaced with the two edges $pv'$ and $v'v$. Cycles without the edge $pv$ remain unchanged.
\medskip

\begin{figure}[H]
\centering
\includegraphics[width=2.15in]{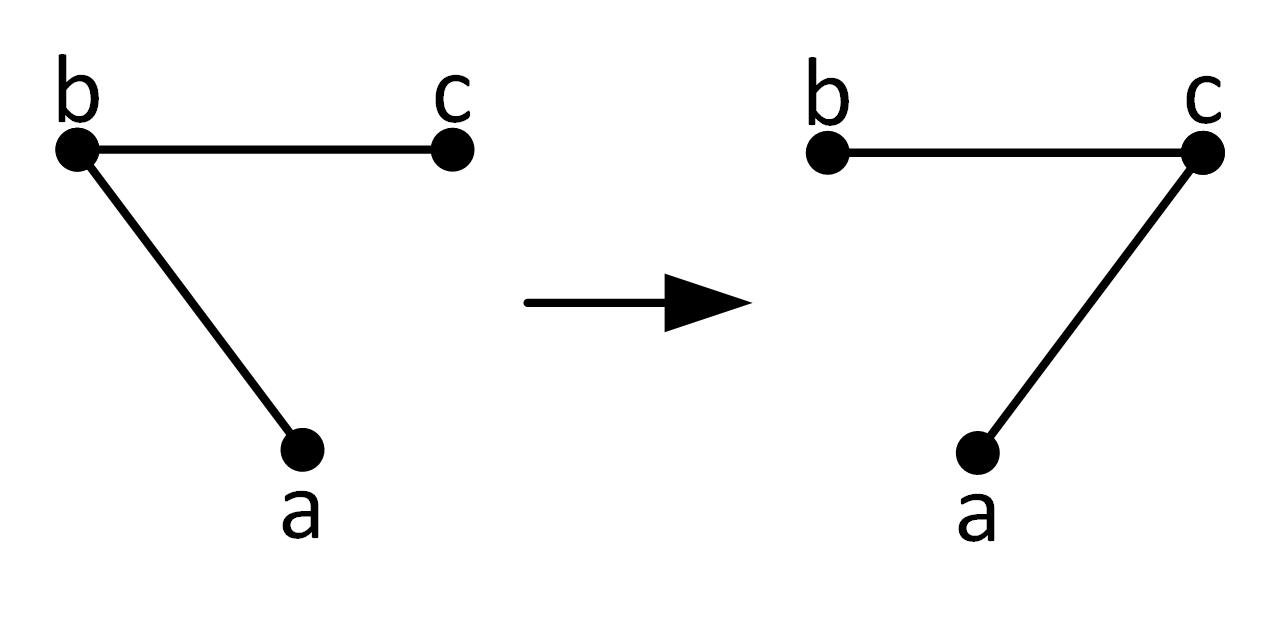}
\caption{ Flipping an edge. \label{flip}}
\end{figure}

The cycles of the graph resulting from step (2) above are more complicated. Suppose $G$ is a graph and consider three vertices $a$, $b$, and $c$ in $G$ where $ab$ and $bc$ are edges, but $ac$ is not an edge. Let $G'$ be the graph formed from $G$ by deleting edge $ab$ and adding edge $ac$. Think of this as ``flipping'' the edge $ab$ to the edge $ac$ as shown in Figure \ref{flip}. Please note that in Figure \ref{vertexsplitsforpropagatingcycles}, this corresponds to removing the edge $qv$ and replacing it with edge $qv'$.

\begin{figure}[H]
\centering
\includegraphics[width=1.5in]{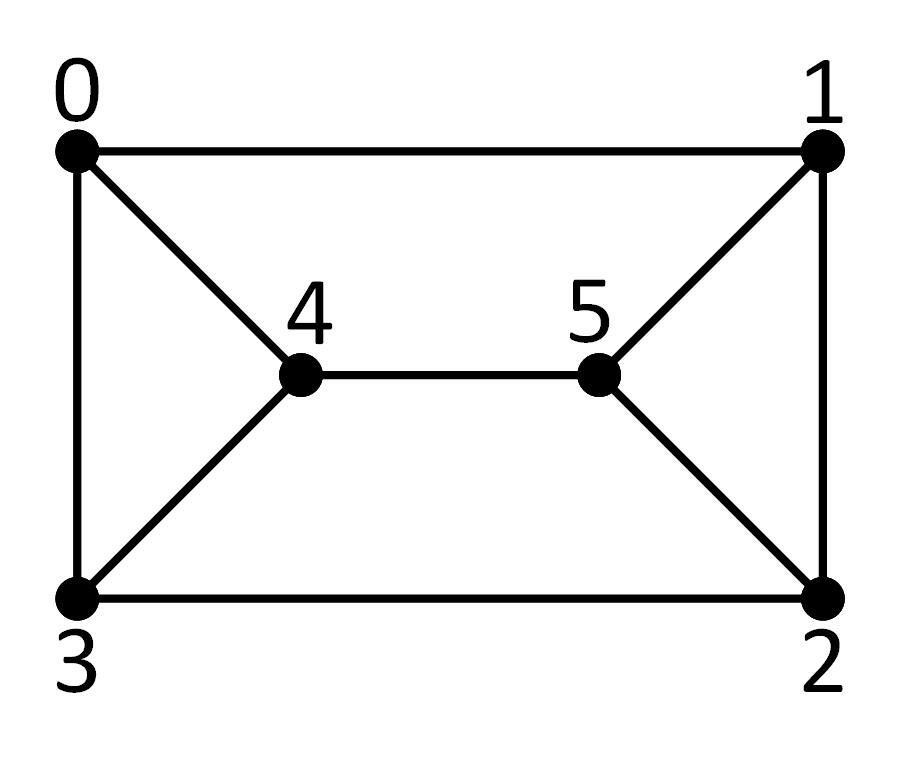}
\caption{The prism graph.\label{prism}}
\end{figure}

Let $C$ be any cycle in $G$ represented by its vertices in order. We may identify cases for determining how individual cycles are changed when $ab$ is replaced with $ac$, by representing a cycle with a ``pattern'' that describes where $a$, $b$,  and $c$ occur in it, if at all. Consider, for example, the cycles of the prism graph with vertices labeled as shown in Figure \ref{prism}:

\begin{equation*}
\begin{split}
&\{{0\, 1\, 5\, 4\, 3\, 0},
{0\, 1\, 2\, 5\, 4\, 3\, 0},
{0\, 1\, 5\, 2\, 3\, 4\, 0},
{0\, 3\, 2\, 1\, 5\, 4\, 0}, \\
&{1\, 2\, 3\, 4\, 5\, 1},
{0\, 1\, 2\, 5\, 4\, 0},
{0\, 1\, 5\, 2\, 3\, 0},
{0\, 1\, 2\, 3\, 4\, 0},
{2\, 3\, 4\, 5\, 2}, \\
&{1\, 2\, 5\, 1},
{0\, 3\, 2\, 5\, 4\, 0},
{0\, 1\, 5\, 4\, 0},
{0\, 3\, 4\, 0},
{0\, 1\, 2\, 3\, 0} \}
\end{split}
\end{equation*}

We identify cycles of the modified graph by following the three steps below, illustrated by the example of the cycle $015430$ taken from the prism graph. Eliminate the redundant final vertex $0$ in the list to obtain $01543$. In this example, let $a=1$, $b=4$, and $c=3$. If none of $a, b, c$ appear in $C$, then there is nothing to do since it remains a cycle in $G'$.

\begin{enumerate}

\item Rotate the list so that $a$ appears first, if it occurs in the cycle, or $b$ if it appears, or $c$ if it appears: $1\, 5\, 4\, 3\, 0$.

\item Replace the vertex numbers associated with $a$, $b$ and $c$ with ``a'', ``b'' and ``c'', respectively: $a\, 5\, b\, c\, 0$.

\item Replace the first sequence of one or more vertices not equal to $a$, $b$ or $c$ with a diamond ($\diamond$), the second if it occurs with a triangle ($\triangle$) and the third, if it occurs, with a square ($\square$): $a \diamond b\, c\, \triangle$.

\end{enumerate}
It helps to think of these steps as symbolic operations:

$15430$

$a5bc0$

$a\diamond  b c \triangle $

\noindent There is no square in the above example. If we start with cycle $012543$  with $a=1$, $b=5$, $c=3$ we get

$125430$

$a2b4c0$

$a\diamond  b \triangle c \square $

\medskip
This procedure will produce different results depending on the orientation used when enumerating the vertices in the cycle; we include all possible patterns in the case-checking in the next result for clarity's sake. Moreover, as explained above, in this representation, $\diamond$, $\triangle$, and $\square$ simply represent sequences of vertices in the cycle other than $a$, $b$, or $c$; the sequences they represent could be of any length. Finally, unlike Lemma \ref{CycleChordingLemma}, there are no connectivity conditions on Lemma \ref{EdgeFlipLemma}.

\begin{lemma}\label{EdgeFlipLemma} {\bf (Edge Flip Lemma)} Let $G$ be a simple graph with $n$ vertices and let $\mathcal C(G)$ be the set of cycles of $G$. Let $a, b, c \in V(G)$ such that $ab, bc \in E(G)$, but $ac \notin E(G)$.
Let $G'$ be the graph obtained from $G$ by replacing $ab$ with a new edge $ac$. The complexity of determining the cycles of $G'$  is $\mathcal O(|\mathcal C(G)|^2n)$.
\end{lemma}

\begin{proof} The cycles of $G'$ can be determined from the cycles of $G$ by analysis of patterns as described~above. First observe that any cycle in $G$ that does not include at least two of the vertices $a$, $b$, and $c$ remains a cycle in $G'$. If a cycle of $G$ does contain at least two of $a$, $b$, and $c$, then we can evaluate how the cycle is affected by the flip from $ab$ to $ac$ based on the cycle's pattern. 

\medskip

We can enumerate all possible patterns by first listing all possible orderings of at least two of $a$, $b$ and $c$: $ab$, $ac$, $bc$, $abc$ and $acb$, and then for each one identifying the possible patterns. Representing cycles in this fashion allows us to distill all of the cycles passing through at least 2 of $a$, $b$ and $c$ in $G$ into 6 cases with a total of 16 subcases for determining how they relate to cycles in $G'$.

\medskip

\noindent {\it Case 1:} $ab$: A pattern containing $a$ and $b$ may or may not include vertices between $a$ and $b$, and may or may not include vertices between $b$ and $a$. This results in four combinations: $ab$, $a\!\diamond \! b$, $a b \diamond$, and~$a\!\diamond \! b \triangle$. Of~these $ab$ is impossible because $G$ has no parallel edges, and therefore a cycle in $G$ must have three~edges. Cycles matching the other three patterns are propagated as follows: \vspace{6pt}

\begin{tabular}{p{3.4in}  p{2in}}
$a\!\diamond \!  b$: If there is a cycle of the form $a\!\diamond \! b$ in $G$ as shown in the left-hand side of the diagram, then when the flip is implemented and $ab$ is replaced with $ac$ in $G'$, $a\!\diamond \! bc$ must be a cycle. In other words $G'$ has a cycle $a\!\diamond \! bc$ in place of cycle $a\!\diamond \! b$. (Cycles in the diagram are indicated with dashed lines.) & \vspace{0pt}\parbox[c]{4em}{\includegraphics[width=1.5in]{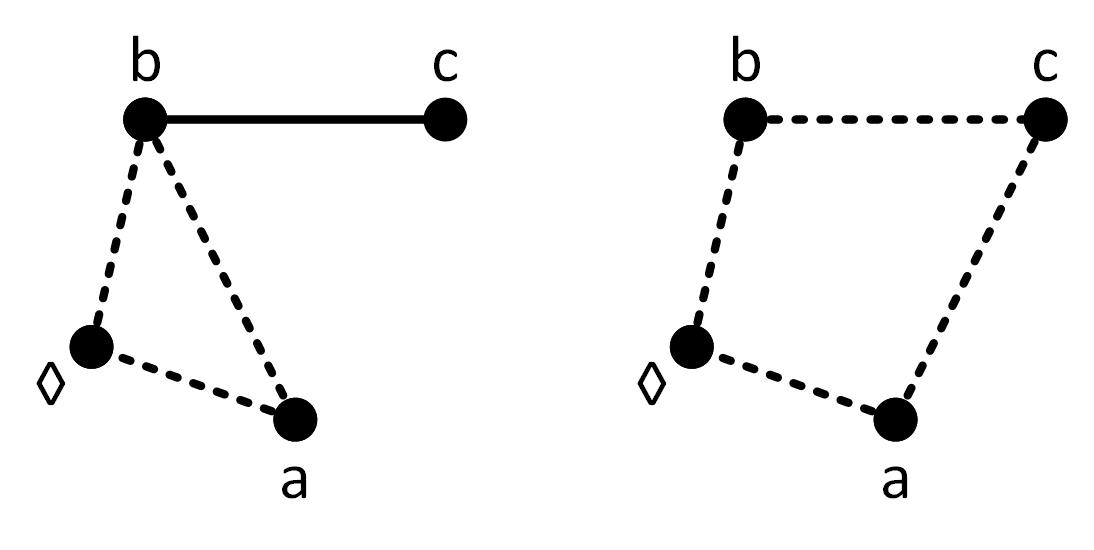}} \\
\end{tabular}

\begin{tabular}{p{3.4in}  p{2in}}
$ab\diamond $: If there is a cycle of the form $ab \diamond$ in $G$, then $G'$ has a cycle $acb \diamond$, which~is $ab \diamond$ with $ab$ replaced with $acb$.& \parbox[c]{4em}{\includegraphics[width=1.5in]{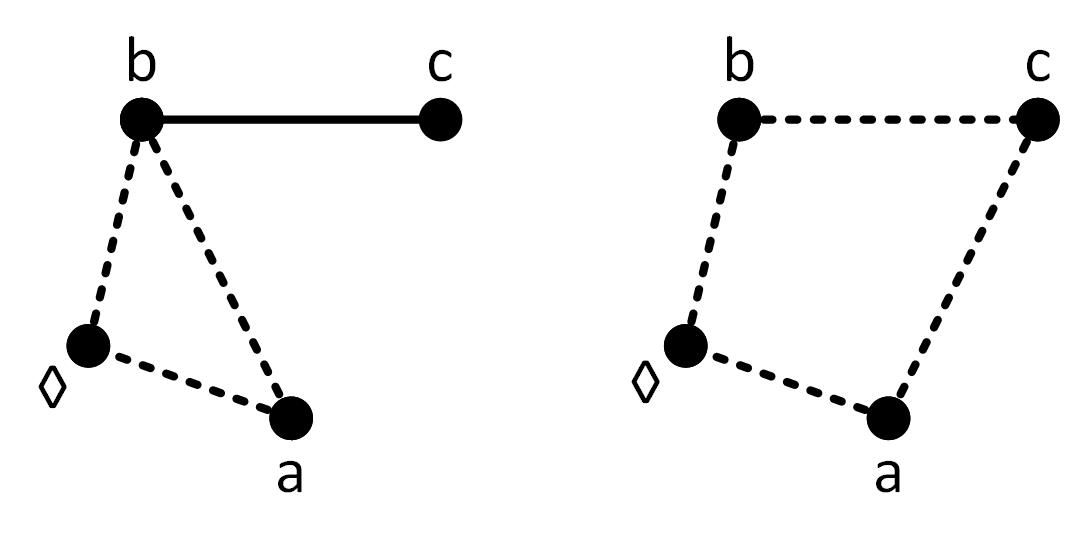}} \\
\end{tabular}

\begin{tabular}{p{3.4in}  p{2in}}
$a\!\diamond \! b\triangle$: This cycle remains a cycle in $G'$. & \parbox[c]{4em}{\includegraphics[width=1.5in]{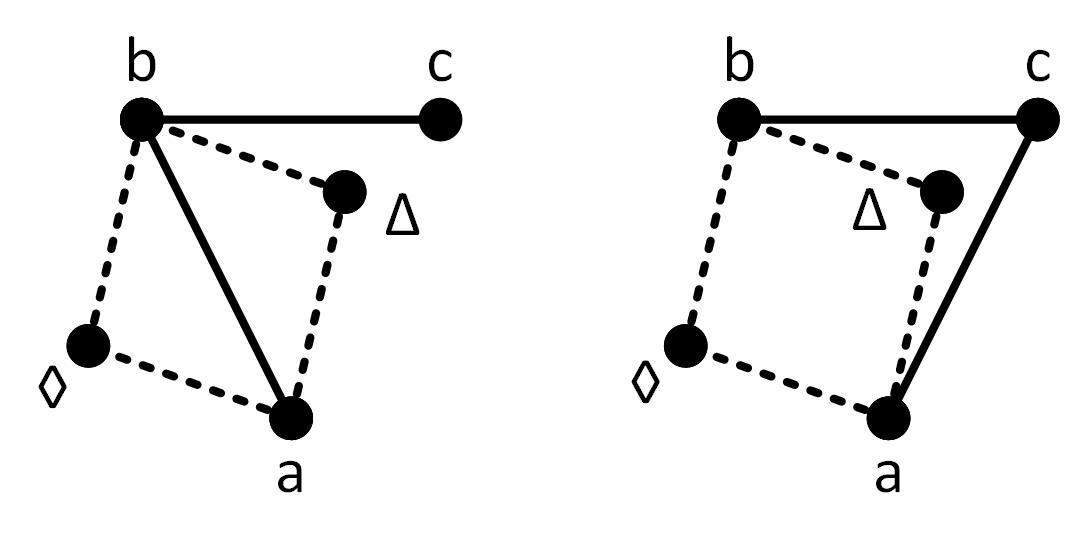}} \\
\end{tabular}

\medskip

\noindent {\it Case 2:} $ac$: The possible patterns containing $a$ and $c$ are $ac$, $a\! \diamond \! c$, $ac \diamond$, and $a\!\diamond \! c \triangle$. In this case, 3 of the 4 patterns are impossible: $ac$ is impossible because $G$ has no parallel edges; $ac \diamond$ and $a\! \diamond \! c$ are impossible because $a$ and $c$ are not adjacent. Cycles matching the remaining pattern are propagated as follows:

\begin{tabular}{p{3.4in}  p{2in}}
$a\!\diamond \! c\triangle $: $G'$ has the same cycle as $G$. Two new cycles emerge also, namely~$a\!\diamond \! c$ and $ac\diamond$, because $ac$ chords the cycle.& \parbox[c]{4em}{\includegraphics[width=1.5in]{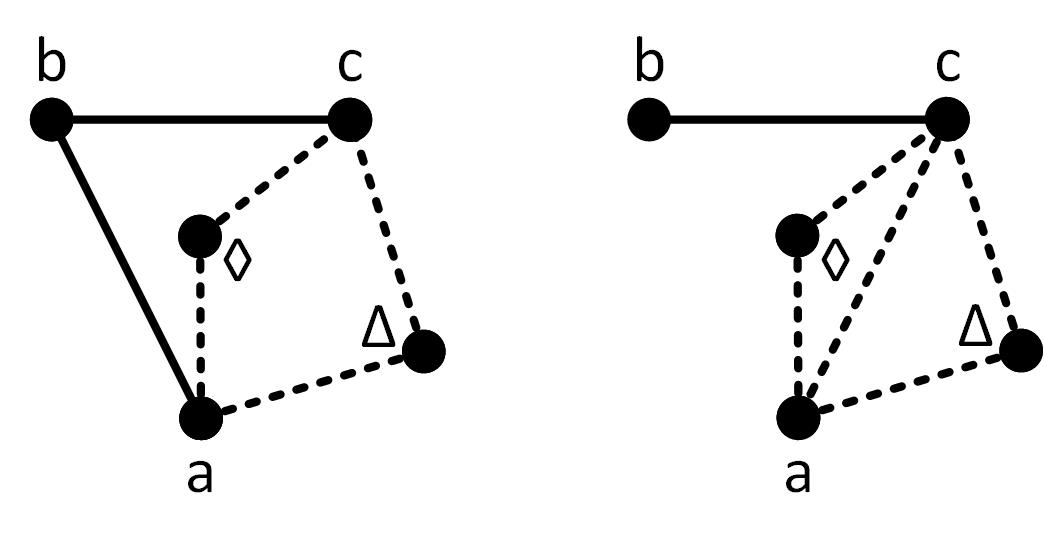}}\vspace{12pt} \\
\end{tabular}

\medskip

\noindent {\it Case 3: } $bc$: The possible patterns containing $b$ and $c$ are $bc$, $b\! \diamond \! c$, $bc \diamond$, and $b\! \diamond \!c \triangle$. In this case, $bc$ is impossible because $G$ has no parallel edges. Cycles matching the other three patterns are propagated with no change:

\begin{tabular}{p{3.4in}  p{2in}}
$b\!\diamond \! c$: This remains a cycle in $G'$.& \parbox[c]{4em}{\includegraphics[width=1.5in]{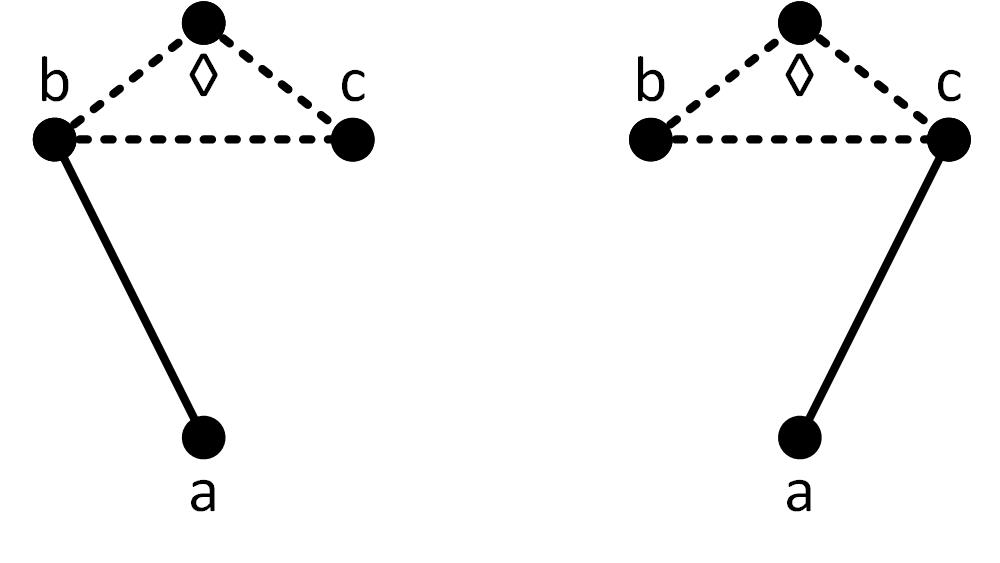}} \\
\end{tabular}

\begin{tabular}{p{3.4in}  p{2in}}
$bc\diamond $: This remains a cycle in $G'$.& \parbox[c]{4em}{\includegraphics[width=1.5in]{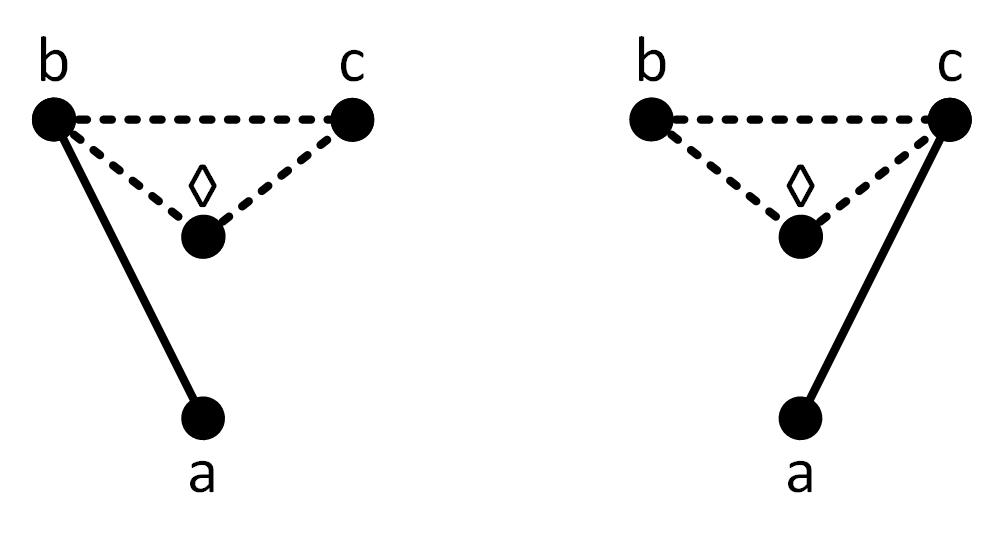}} \\
\end{tabular}

\begin{tabular}{p{3.4in}  p{2in}}
$b\!\diamond \! c\triangle$: This remains a cycle in $G'$.& \parbox[c]{4em}{\includegraphics[width=1.5in]{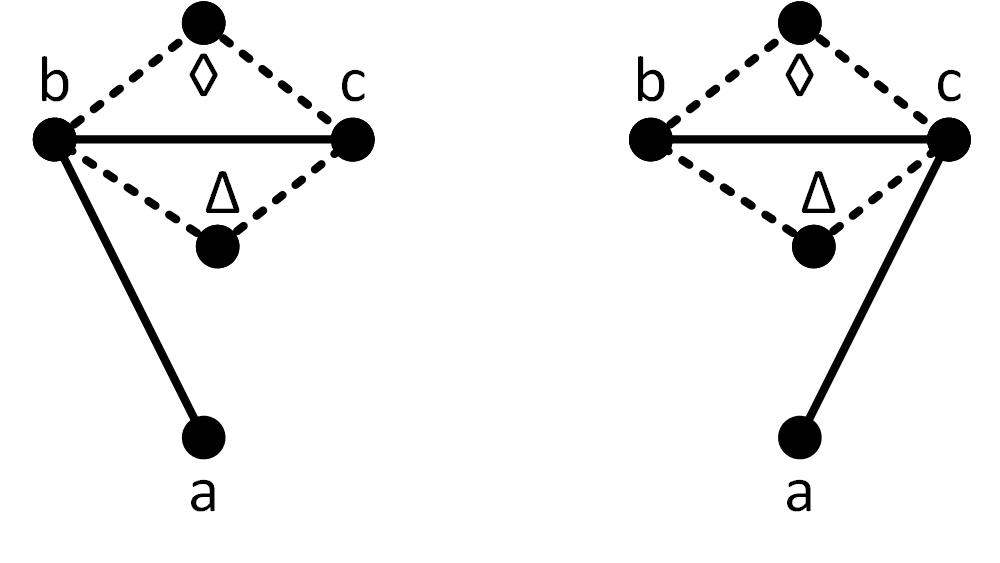}} \\
\end{tabular}

\medskip

\noindent {\it Case 4:} $abc$: The eight possible patterns containing $a$, $b$, and $c$ in order are $abc$, $a\! \diamond \! bc$, $ab\! \diamond \!c$, $a\! \diamond \!b \triangle c$, $abc \diamond$, $a\! \diamond \! bc \triangle$, $ab\! \diamond \! c \triangle$, and $a\! \diamond \! b \triangle c \square$. In this case, four patterns, $abc$, $a\! \diamond \! bc$, $ab\! \diamond \!c$, and $a\! \diamond \!b \triangle c$ are all impossible because $a$ and $c$ are not adjacent in $G$. Cycles matching the other four patterns are propagated as~follows:

\begin{tabular}{p{3.4in}  p{2in}}
$abc\diamond $: If $G$ has a cycle of the form $abc \diamond$, then $G'$ has a cycle $ac \diamond$, which is $abc \diamond$ with $abc$ replaced with $ac$.& \parbox[c]{4em}{\includegraphics[width=1.5in]{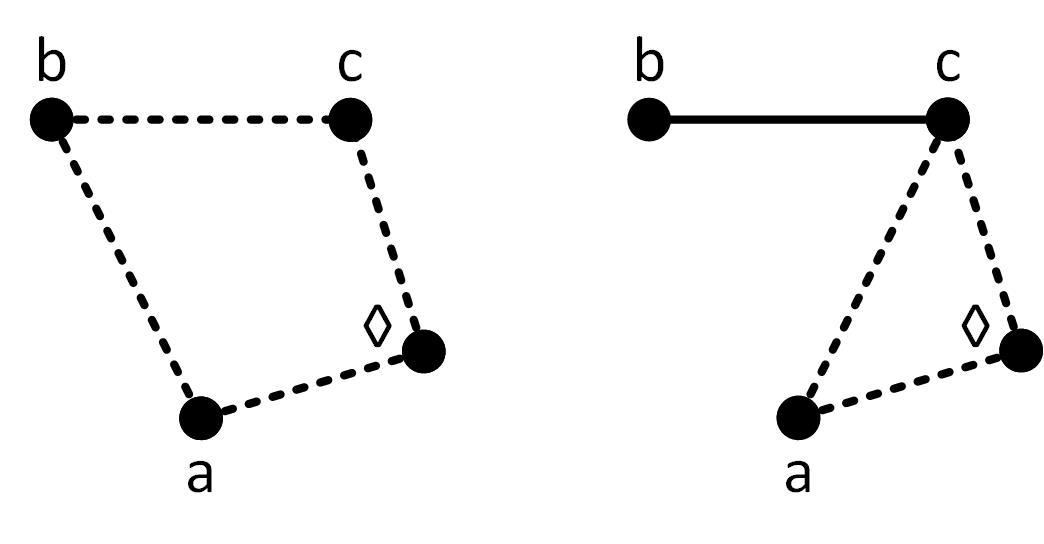}} \\
\end{tabular}

\begin{tabular}{p{3.4in}  p{2in}}
$a\!\diamond \! bc\triangle $: This remains a cycle in $G'$.& \parbox[c]{4em}{\includegraphics[width=1.5in]{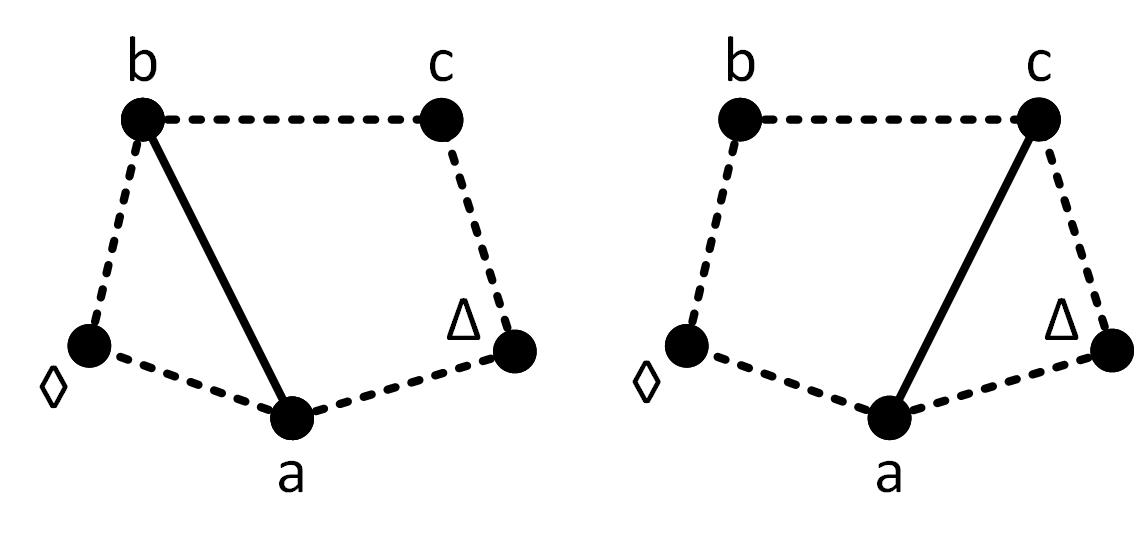}} \\
\end{tabular}

\begin{tabular}{p{3.4in}  p{2in}}
$ab\!\diamond \! c\triangle $: If $G$ has a cycle of the form $ab \! \diamond \! c \triangle$, then it will be replaced in $G'$ with two cycles: $b\!\diamond \! c$ and $ac\triangle$.& \parbox[c]{4em}{\includegraphics[width=1.5in]{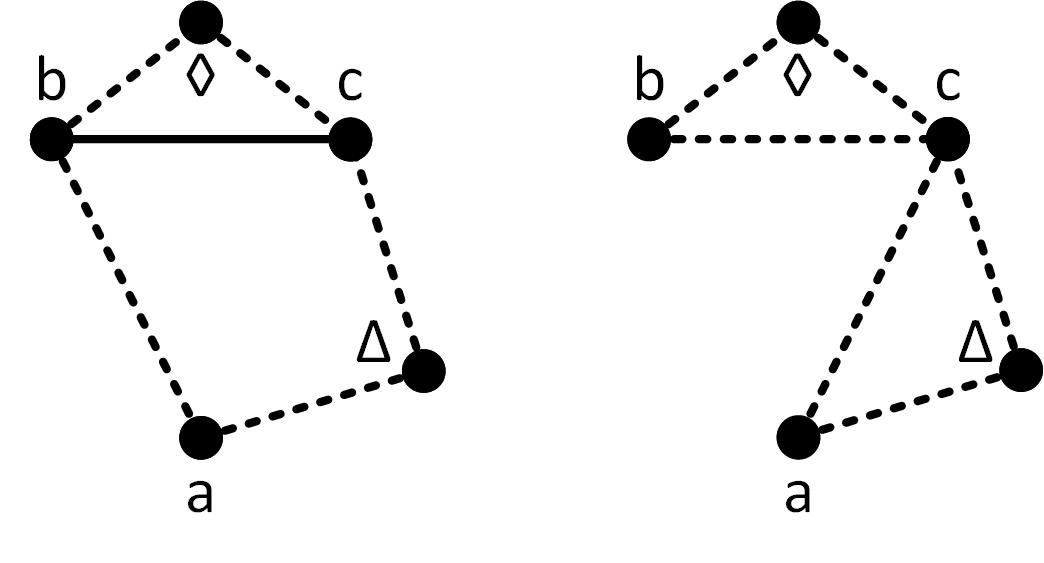}} \\
\end{tabular}

\begin{tabular}{p{3.4in}  p{2in}}
$a\!\diamond \! b\triangle c\square $: This remains a cycle in $G'$.& \parbox[c]{4em}{\includegraphics[width=1.5in]{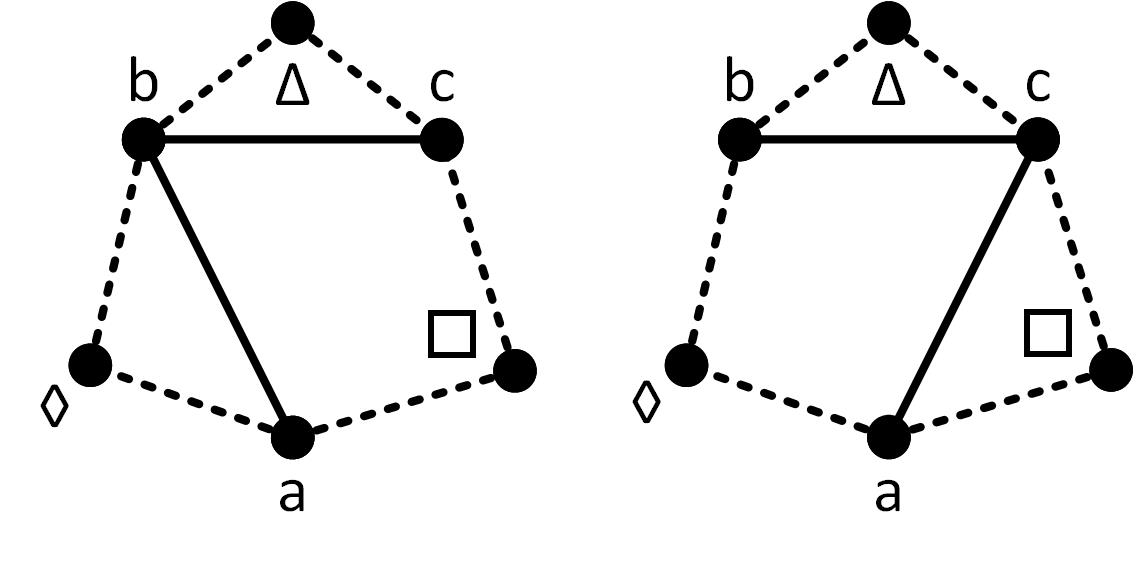}} \\
\end{tabular}

\medskip

\noindent {\it Case 5:} $acb$: The eight possible patterns containing $a$, $c$, and $b$ in order are $acb$, $a\! \diamond \! cb$, $ac\! \diamond \! b$, $a\! \diamond \! c \triangle  b$, $acb \diamond$, $a\! \diamond \! cb \triangle$, $ac\! \diamond \!b \triangle$, and $a\! \diamond \! c \triangle b \square$. In this case, four patterns, $acb$, $ac\! \diamond \! b$, $acb \diamond$, and $ac\! \diamond \!b \triangle$ are impossible because $a$ and $c$ are not adjacent in $G$. Cycles matching the other four patterns are propagated as~follows:

\begin{tabular}{p{3.4in}  p{2in}}
$a\!\diamond \! cb$: If $G$ has a cycle of the form $a\! \diamond \! cb$, then $G'$ will have a cycle of the form $a\! \diamond \! c$, which is the original cycle with $cba$ replaced with $ca$.& \parbox[c]{4em}{\includegraphics[width=1.5in]{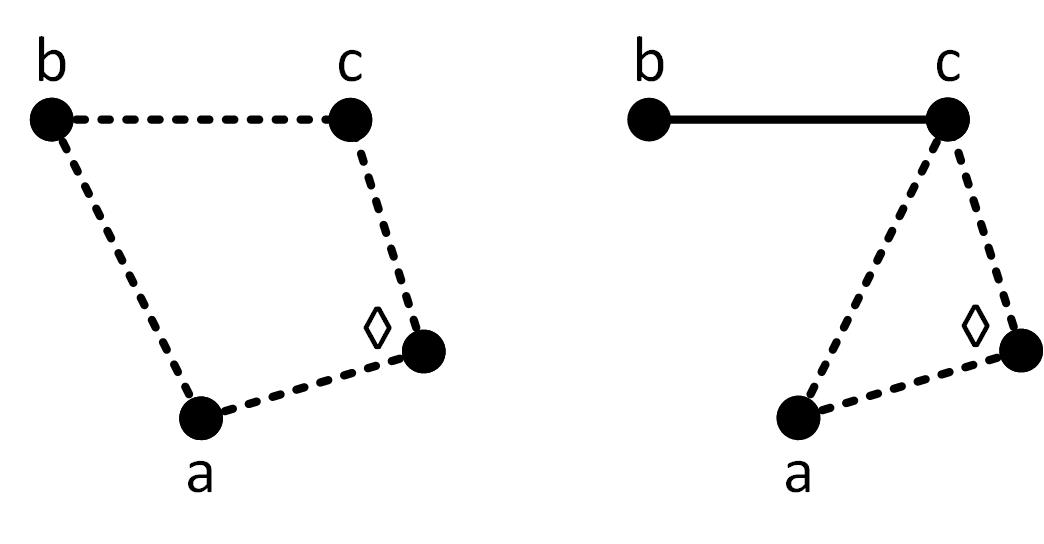}} \\
\end{tabular}

\begin{tabular}{p{3.4in}  p{2in}}
$a\!\diamond \! c\triangle b $: If $G$ has a cycle of the form $a\! \diamond \! c\triangle b$, then $G'$ will have cycles of the form $b \triangle c$ and $a \diamond c$ in its place.& \parbox[c]{4em}{\includegraphics[width=1.5in]{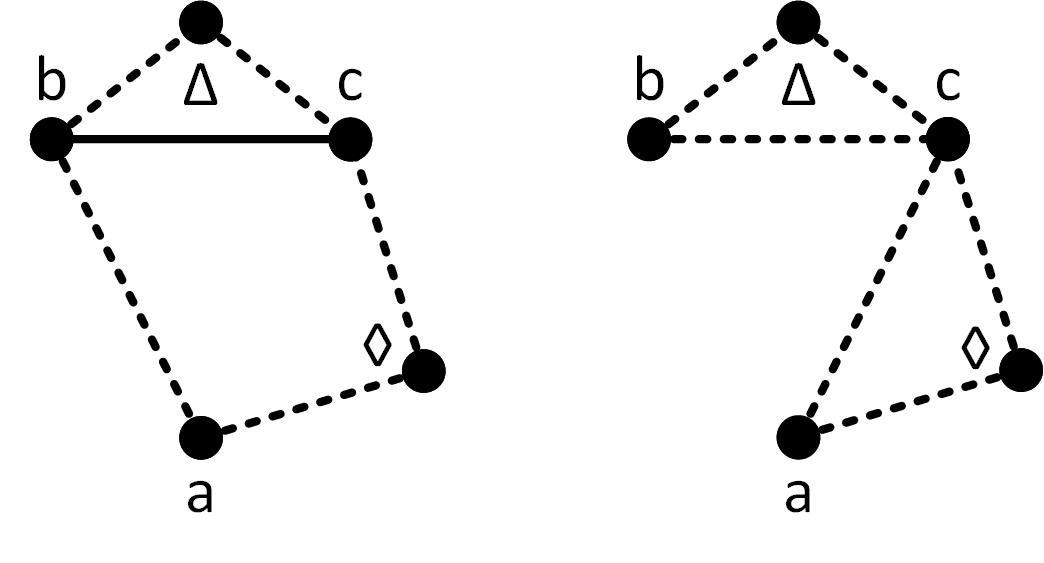}} \\
\end{tabular}

\begin{tabular}{p{3.4in}  p{2in}}
$a\!\diamond \! cb\triangle $: This remains a cycle in $G'$.& \parbox[c]{4em}{\includegraphics[width=1.5in]{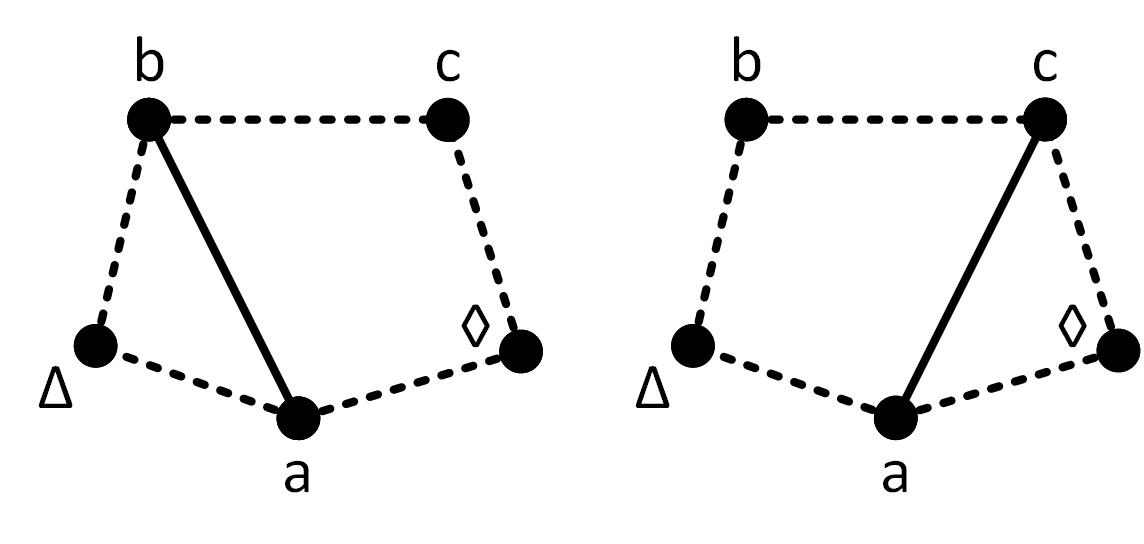}} \\
\end{tabular}

\begin{tabular}{p{3.4in}  p{2in}}
$a\!\diamond \! c\triangle b\square $: This remains a cycle in $G'$.& \parbox[c]{4em}{\includegraphics[width=1.5in]{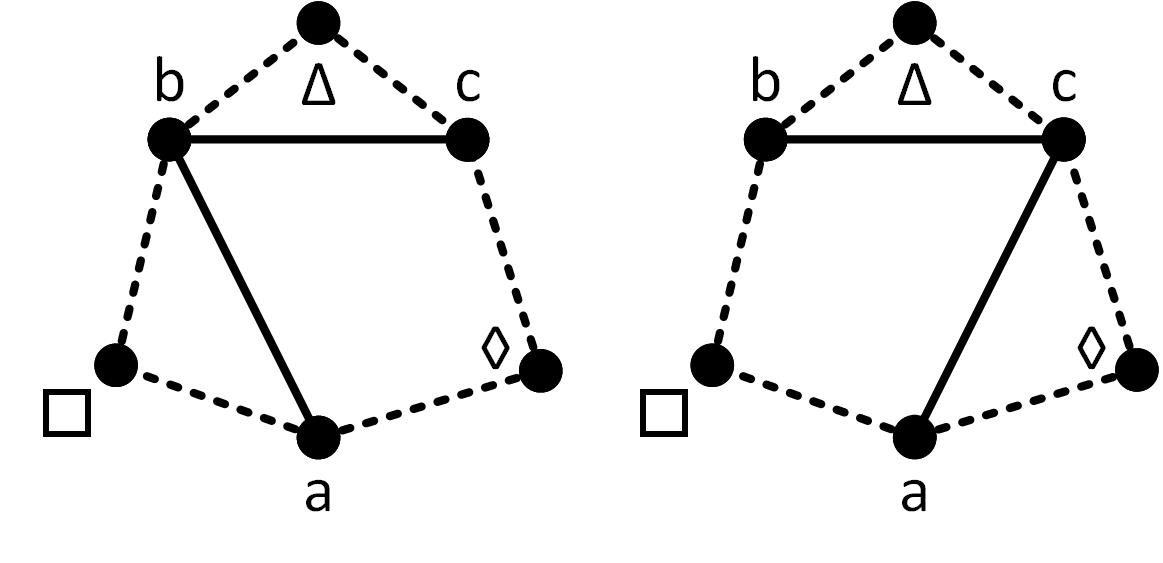}} \\
\end{tabular}

\noindent {\it Case 6:} There is one additional case in which two cycles in $G$ result in one cycle in $G'$ after the flip~operation:  \vspace{6pt}

\begin{tabular}{p{3.4in}  p{2in}}
Two cycles in $G$ which share the common vertex $b$, share no other common vertices and for which the edge $ab$ lies in one cycle and the edge $bc$ lies in the other; that is a pair of cycles with patterns $ab\diamond$ and $b \triangle c$, correspond to one cycle in $G'$ of the form $a\!\diamond \! b\triangle c$.& \parbox[c]{4em}{\includegraphics[width=1.5in]{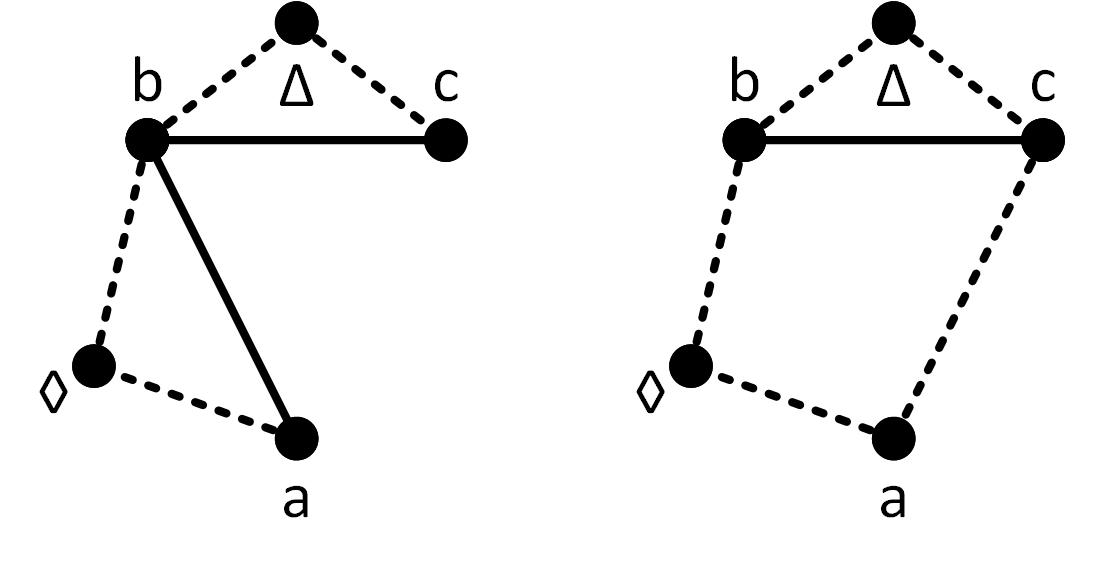}}\\
\end{tabular}

\medskip
\vspace{6pt}
In all but the last case, an existing cycle has to be traversed to produce a new cycle making it an $\mathcal O(n)$ operation because a cycle may contain at most $n$ vertices. Without the last case, because each cycle has to be traversed the complexity would be $\mathcal O(|\mathcal C(G)|n)$. The last case requires consideration of every pair of cycles which is $\mathcal O(|\mathcal C(G)|^2n)$. \end{proof}

\medskip

Now, using Lemmas \ref{CycleChordingLemma} and   \ref{EdgeFlipLemma} we can establish bounds on the complexity of identifying the cycles of a graph obtained by one of operations D1, D2, and D3, in terms of the cycles of the original graph.

\medskip

\begin{theorem} \label{cycle-prop-theorem} Let $G'$ be a simple graph obtained from a smaller $3$-connected graph $G$ by one of operations D1, D2, and D3. Let $n$ be the number of vertices in $G$ and let $c$ be the number of cycles of $G$. Then the cycles of $G'$ can be obtained from the cycles of $G$ by a method with complexity $\mathcal O(c^2 n)$.
\end{theorem}

\begin{proof} Using Theorem \ref{maintheorem-1}, operation D1 can be expressed as an edge addition, followed by an edge subdivision, followed by an edge flip. By Lemmas \ref{CycleChordingLemma} and   \ref{EdgeFlipLemma}, the complexities for these individual steps are $\mathcal O(c n)$, $\mathcal O(c)$, and $\mathcal O(c^2 n)$, respectively, so the overall complexity is $\mathcal O(c^2 n)$. Similarly, operation D2 can be expressed as an edge addition, followed by two edge subdivisions and edge flips, and operation D3 can be expressed as two edge additions followed by an edge subdivision and an edge flip, so the overall complexity of propagating the list of cycles for D2 and D3 is also $\mathcal O(c^2 n)$.
\end{proof}

\bigskip

%%%%%%%%%%%%%%%%%%%%%%%%%%%%%%%

\section{\bf Cycle Propagation Algorithm}\label{cycle-propagation-section}

\bigskip

To check when operations D1, D2, and D3 can be performed, we must check for chording paths between specific vertices. Let $G$ be a graph with $n$ vertices and let $\mathcal C$ denote its set of cycles.  Since~we will be considering only one input graph at a time, in this section $\mathcal C$ is unambiguous.  Given the set of cycles of the input graph $G$, Lemmas \ref{CycleChordingLemma} and \ref{EdgeFlipLemma} allow us to define three procedures \textproc{ApplyFlipEdge},  \textproc{ApplyAddEdge}, and \textproc{ApplySubdivideEdge}, which compute the set of cycles of a graph modified by flipping an edge, adding an edge, or subdividing an edge. This allows us to propagate the set of cycles for the new graphs we produce in terms of the cycles of $G$, which in turn allows us to avoid the NP-complete problem of computing all the cycles of a graph.

\medskip

The procedures described in this section rely on two simpler procedures, \textproc{ChordCycle} and \textproc{Chords}, which are defined informally below. 

\medskip

\begin{itemize}
\item \textproc{ChordCycle}($C, v_1, v_2$): This is the procedure described in Lemma \ref{CycleChordingLemma}. Given a cycle $C$ and a pair of vertices $v_1$ and $v_2$ which both occur in the cycle, but are not adjacent, return a pair of cycles generated as follows: The first cycle is created by starting at $v_1$, then proceeding to $v_2$, and then following the existing cycle from $v_2$ back to $v_1$. The second cycle begins with $v_1$ and proceeds along the existing cycle until $v_2$ is reached, and then returns to $v_1$. That is, if the existing cycle consists of vertices $$w_1, w_2, \ldots, w_i, v_1, w_{i + 1}, \ldots, w_j, v_2, w_{j+1}, \ldots, w_k,$$ then the two new cycles will  consist of vertices $$v_2, w_{j + 1}, \ldots, w_k, w_1, \ldots, v_1$$ and $$v_1, w_{i + 1}, \ldots, w_j, v_2.$$ Its complexity is $\mathcal O(n)$.

\medskip

\item \textproc{Chords}($C, v_1, v_2$): This procedure simply determines whether an edge chords a cycle. Given~a~cycle $C$ and two vertices $v_1$ and $v_2$, determine whether there is an edge $v_1 v_2$ that chords $C$. Its complexity is $\mathcal O(n)$.
\end{itemize}

\medskip

Next, we present the three procedures \textproc{ApplyFlipEdge}, \textproc{ApplyAddEdge} and \textproc{ApplySubdivideEdge}, respectively. 

\medskip

\begin{enumerate}

\item \textproc{ApplyFlipEdge}($\mathcal C, a, b, c$): This procedure is also described informally, as its steps closely follow the cases in the proof of Lemma \ref{EdgeFlipLemma}.  Given the set $\mathcal C$ of cycles of a graph $G$ and vertices $a$, $b$, and $c$ with edges $ab$ and $bc$, but no edge $ac$, \textproc{ApplyFlipEdge} generates the list of cycles of the graph $G'$ produced by replacing edge $ab$ with a new edge $ac$, using the procedure outlined with Lemma \ref{EdgeFlipLemma}. Its complexity is $\mathcal O(|\mathcal C|^2n)$.

\medskip

\item \textproc{ApplyAddEdge}($\mathcal C, a, b$): This procedure computes the resulting cycles when adding an edge to a graph, where $\mathcal C$ is the set of cycles and $a$ and $b$ are two non-adjacent vertices. It creates the new set of cycles by adding the cycles from the old graph and the new cycles created by \textproc{ChordCycle} for any cycle that is chorded by $ab$. Its complexity is $\mathcal{O}(|\mathcal C|n)$. Pseudocode is shown in Algorithm \ref{alg01}.  

\medskip
\begin{algorithm}
\caption{Compute cycles for an edge addition}\label{alg01}
{\scriptsize
\begin{algorithmic}[1]
\Procedure{ApplyAddEdge}{$\mathcal C, a, b$}\vspace{3pt}
\State $S \gets \phi$\vspace{3pt}
\For{$C \in \mathcal C$}\vspace{3pt}
	\State{$S \gets S \cup \{C\}$}\vspace{3pt}
	\If{$\Call{Chords}{C, a, b}$}\vspace{3pt}
		\State{$C_1, C_2 \gets \Call{ChordCycle}{C, a, b}$}\vspace{3pt}
		\State{$S \gets S \cup \{C_1, C_2\}$}\vspace{3pt}
	\EndIf\vspace{3pt}
\EndFor\vspace{3pt}
\EndProcedure\vspace{3pt}
\end{algorithmic}
}
\end{algorithm}

\medskip

\item \textproc{ApplySubdivideEdge}($\mathcal C, a, b, c$): This procedure computes cycles resulting when subdividing an~edge, where $\mathcal C$ is the set of cycles of the graph, $a$ and $b$ are two adjacent vertices, and vertex $c$ is being added to subdivide the edge $ab$. Its complexity is $\mathcal{O}(|\mathcal C|n)$. Pseudocode is shown in Algorithm \ref{alg02}.  

\medskip

\begin{algorithm}
\caption{Compute cycles for subdividing an edge}\label{alg02}
{\scriptsize
\begin{algorithmic}[1]
\Procedure{ApplySubdivideEdge}{$\mathcal C$, $a$, $b$, $c$} \vspace{3pt}
\State $S \gets \phi$\vspace{3pt}
\For {$C \in \mathcal C$}\vspace{3pt}
	\If{$C$ contains edge $(ab)$}\vspace{3pt}
		\State Add cycle $C$ with $a, b$ replaced with $a, c, b$ to $S$\vspace{3pt}
	\Else\vspace{3pt}
		\State Add $C$ to $S$\vspace{3pt}
	\EndIf\vspace{3pt}
\EndFor\vspace{3pt}
\EndProcedure\vspace{3pt}
\end{algorithmic}
}
\end{algorithm}
\end{enumerate}

%%%%%%%%%%%%%%%%%%%%%%%%%%%%%%%

\section{Isomorph-Free Graph Construction}\label{InfiniteBookshelf-section} 
\bigskip

This section is further broken into three subsections.  

\subsection{ {Organizing Graph Construction to Minimize Isomorphism Checking}} 

When we apply operation D1 to a graph, we end up with a graph that has two more edges and one more vertex. When we apply operation D2 to a graph, we end up with a graph that has three more edges and two more vertices. When we apply operation D3 to a graph, we end up with a graph that has three more edges and one more vertex. This creates a problem if we want to avoid generating isomorphic graphs, because we have to keep track of graphs of different sizes at the same~time. To~prevent this, we want to focus on doing everything we need to do with graphs with one particular number of edges and vertices all at once. 
In particular, if we consider operations D1, D2, and D3 as~algorithms, then:

\begin{itemize}

\item D1 takes a graph $G$ with $n$ vertices and $m$ edges, a vertex $x \in V(G)$ and an edge $ab \in E(G)$ as input, and produces a graph $G' = (G + e) \circ f_c$ with $n+1$ vertices and $m+2$ edges (see~Theorem~\ref{maintheorem-1}~(i));

\item D2 takes a graph $G$ with $n$ vertices and $m$ edges, and two edges $ab, cd \in E(G)$ as input, and produces a graph $G' = (G + e) \circ \{f_c, f_b\}$ with $n+2$ vertices and $m+3$ edges (see Theorem \ref{maintheorem-1} (ii)); and

\item D3 takes a graph $G$ with $n$ vertices and $m$ edges, and three vertices $x, y, z \in V(G)$ as input, and~produces a graph $G' = (G + \{e_1, e_2\}) \circ f$ with $n+1$ vertices and $m+3$ edges (see~Theorem~\ref{maintheorem-1}~(iii)).
\end{itemize}

\medskip

\begin{figure}[H]
\centering
\includegraphics[width=5.0in]{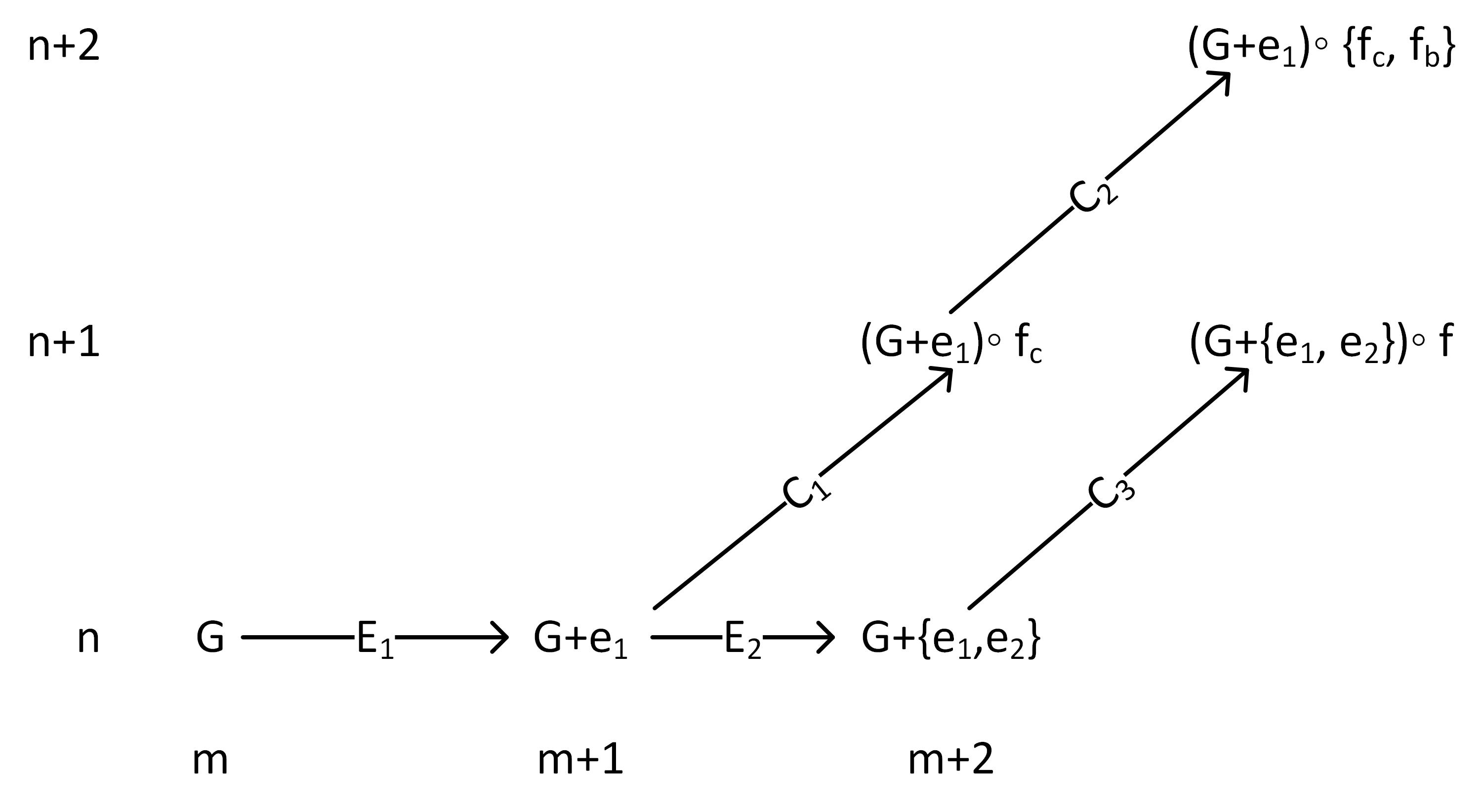}
\caption{Minimally $3$-connected graph construction.\label{flow}}
\end{figure}

Figure \ref{flow} outlines the process of applying operations D1, D2, and D3 to an individual graph. The~specific procedures \textproc{E1}, \textproc{E2}, \textproc{C1}, \textproc{C2}, and \textproc{C3} will be detailed in Section \ref{graph-construction-section}.

To avoid generating graphs that are isomorphic to each other, we wish to maintain a~list of generated graphs and check newly generated graphs against the list to eliminate those for which isomorphic duplicates have already been generated. We immediately encounter two problems with this approach: checking whether a pair of graphs is isomorphic is a computationally expensive operation; and the number of graphs to check grows very quickly as the size of the graphs, both in terms of vertices and edges, increases.
\medskip

The first problem can be mitigated by using McKay's nauty system  [\ref{McKay1981}] ({available for download at \texttt{http://pallini.di.uniroma1.it/}}) to generate certificates for each graph. The nauty certificate function $\phi$  produces a data artifact from a graph in such a way that $\phi(G_1) = \phi(G_2)$ if and only if $G_1 \cong G_2$. Thus~we can reduce the problem of checking isomorphism to the problem of generating certificates, and then compare a newly generated graph's certificate to the set of certificates of graphs already~generated. 

\medskip

\begin{figure}[H]
\centering
\includegraphics[width=3.15in]{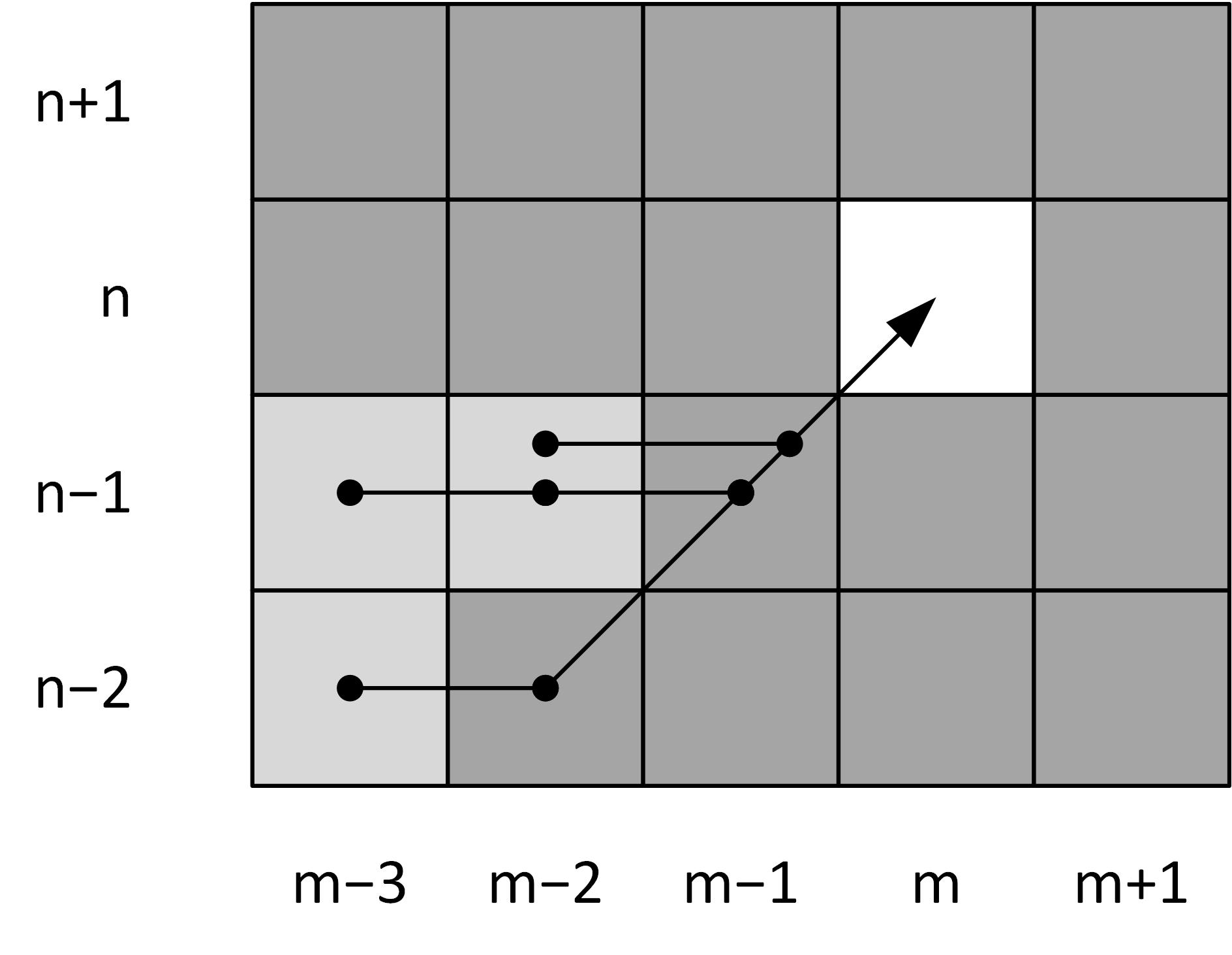}
\caption{Infinite bookshelf approach to minimizing isomorphism checking.\label{infinite-bookshelf}}  
\end{figure}

The second problem can be mitigated by a change in perspective. While Figure \ref{flow} demonstrates how a single graph will be treated by our process, consider Figure \ref{infinite-bookshelf}, which we refer to as the ``infinite bookshelf''. As the entire process of generating minimally 3-connected graphs using operations D1, D2, and D3 proceeds, with each operation divided into individual steps as described in Theorem \ref{maintheorem-1}, the set of all generated graphs with $n$ vertices and $m$ edges will contain both ``finished'', minimally $3$-connected graphs, and ``intermediate'' graphs generated as part of the process. What does this set of graphs look like?

\begin{itemize}

\item For operation D1, the set may include graphs of the form $G + e_1$, where $G$ has $n$ vertices and $m-1$ edges, and graphs of the form $(G + e_1) \circ f$, where $G$ has $n-1$ vertices and $m-2$ edges.

\item For operation D2, the set may include graphs of the form $G + e_1$, where $G$ has $n$ vertices and $m-1$ edges, graphs of the form $(G + e_1) \circ f$, where $G$ has $n-1$ vertices and $m-2$ edges, and~graphs of the form $(G + e_1) \circ \{f_c, f_b\}$, where $G$ has $n-2$ vertices and $m-3$ edges.

\item For operation D3, the set may include graphs of the form $G + e_1$ where $G$ has $n$ vertices and $m-1$ edges, graphs of the form $G + \{e_1, e_2\}$, where $G$ has $n$ vertices and $m-2$ edges, and graphs of the form $(G + \{e_1, e_2\}) \circ f$, where $G$ has $n-1$ vertices and $m-3$ edges.

\end{itemize}

When generating graphs, by storing some data along with each graph indicating the steps used to generate it, and by organizing graphs into subsets, we can generate all of the graphs needed for the algorithm with $n$ vertices and $m$ edges in one batch. Organized in this way, we only need to maintain a~list of certificates for the graphs generated for one ``shelf'', and this list can be discarded as soon as processing for that shelf is complete. We do not need to keep track of certificates for more than one shelf at a time.

\medskip

Section \ref{isomorph-free-subsection} breaks down the graphs in one shelf formally by their place in operations D1, D2, and~D3. Section \ref{graph-construction-section} then describes how the procedures for each shelf work and interoperate.

\bigskip

\subsection{ The Algorithm Is Isomorph-Free}\label{isomorph-free-subsection} 

To make the process of eliminating isomorphic graphs by generating and checking nauty certificates more efficient, we organize the operations in such a way as to be able to work with all graphs with a fixed vertex count $n$ and edge count $m$ in one batch. Specifically, for an $m, n$ combination, we define sets $A^{(*)}_{m, n}$, where $*$ represents 0, 1, 2, or 3, $B_{m, n}$ and $C_{m, n}$ as follows:

\medskip

\begin{itemize}
\item $A^{(0)}_{m, n}$ only ever contains of the ``root'' graph; i.e., the prism graph. So for values of $m$ and $n$ other than $9$ and $6$, $A^{(0)}_{m, n} = \varnothing$. All graphs in $A^{(0)}_{m, n}$, $A^{(1)}_{m, n}$, $A^{(2)}_{m, n}$, and $A^{(3)}_{m, n}$ are minimally 3-connected.

\medskip

\item $B_{m, n}$ consists of graphs generated by adding an edge to a minimally 3-connected graph with $m-1$ vertices and $n$ edges.

\item $C_{m, n}$ consists of graphs generated by adding an edge to a graph in $B_{m-1, n}$ that is incident with the edge added to form the input graph.

\item $A^{(1)}_{m, n}$ consists of graphs generated by splitting a vertex in a graph in $B_{m-1, n-1}$ that is incident to the edge added to form the input graph, after checking for 3-compatibility.

\item $A^{(2)}_{m, n}$ consists of graphs generated by splitting a vertex in a graph in $A^{(1)}_{m-1, n-1}$ that is incident to the same edge as the vertex split to form the input graph, after checking for 3-compatibility.

\item $A^{(3)}_{m, n}$ consists of graphs generated by splitting a vertex in a graph in $C_{m-1, n-1}$ that is incident to the two edges added to form the input graph, after checking for 3-compatibility.

\end{itemize}

Then, beginning with $m = 10$ and $n = 6$, we construct graphs in $C_{m,n}$, $B_{m,n}$, $A^{(1)}_{m, n}$, $A^{(2)}_{m,n}$ and $A^{(3)}_{m,n}$, in that order, from input graphs with $m-1$ vertices and $n$ edges, and with $m-1$ vertices and $n-1$~edges. As graphs are generated in each step, their certificates are also generated and stored. Any new graph with a certificate matching another graph already generated, regardless of the step, is~discarded, so that the full set of generated graphs is pairwise non-isomorphic. At the end of processing for one value of $n$ and $m$ the list of certificates is discarded.

\medskip

\begin{figure}[H]
\centering
\includegraphics[width=3.9in]{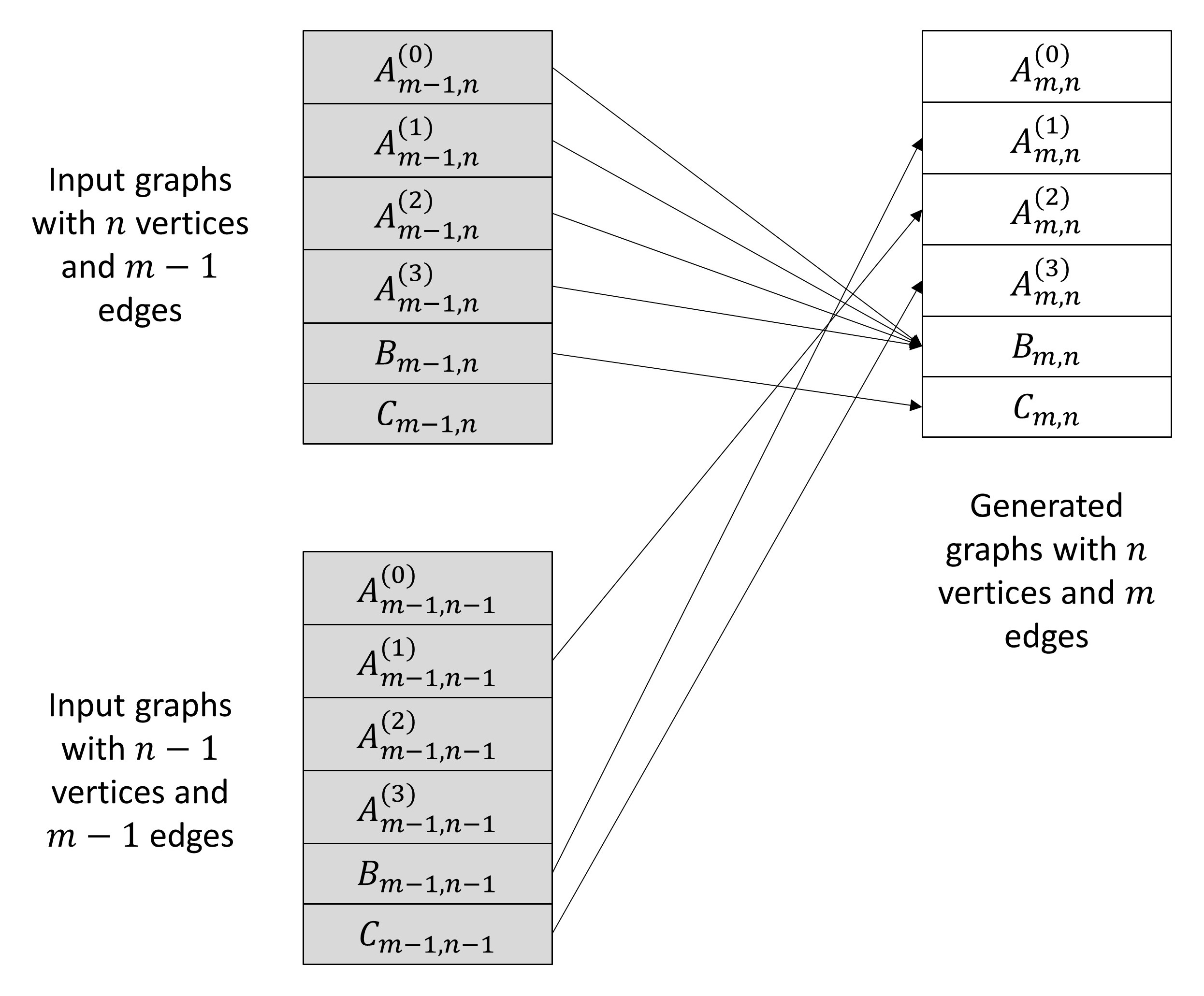}
\caption{Graph generation steps. \label{bookshelf}}
\end{figure}

For any value of $n$, we can start with $m = n + 4$ and proceed until no more graphs or generated or, when $n \ge 8$, when $m = 3n - 7$. By Theorem \ref{HalinTheorem}, no further minimally 3-connected graphs will be found after $m = 3n - 9$ when $n \ge 8$; however we still need to generate single- and double-edge additions to be used when considering graphs with $n+1$ vertices. Proceeding in this fashion, at any time we only need to maintain a list of certificates for the graphs for one value of $m$ and $n$. The generation sources and targets are summarized in Figure \ref{bookshelf}, which shows how the graphs with $n$ vertices and $m$ edges, in~the upper right-hand box, are generated from graphs with $n$ vertices and $m-1$ edges in the upper left-hand box, and graphs with $n-1$ vertices and $m-1$ edges in the lower left-hand box.

 \bigskip
%%%%%%%%%%%%%%%%%%%%%%%%%%%%%%%

\subsection{Infinite Bookshelf Algorithm}\label{graph-construction-section}
This subsection contains a detailed description of the algorithms used to generate graphs, implementing the process described in Section \ref{isomorph-free-subsection}.  

\medskip

 The process of computing $(G + e) \circ f_c$, $(G + e) \circ \{f_c, f_b\}$, and $(G + \{e_1, e_2\}) \circ f$ is broken down into individual procedures \textproc{E1}, \textproc{E2}, \textproc{C1}, \textproc{C2}, and \textproc{C3}, each of which operates on an input graph with one less edge, or one less edge and one less vertex, than the graphs it produces.
The procedures are implemented using the following component steps, as illustrated in Figure \ref{flow}:

\begin{itemize}

\item Procedure \textproc{E1} is applied to graphs in $A^{(*)}_{m-1, n}$, which are minimally 3-connected, to generate all possible single edge additions $G+e$ given an input graph $G$. This is the first step for operations D1, D2, and D3, as expressed in Theorem \ref{maintheorem-1}. The cycles of the output graphs are constructed from the cycles of the input graph $G$ (which are carried forward from earlier computations) using \textproc{ApplyAddEdge}. The results, after checking certificates, are added to $B_{m, n}$.
\medskip

\item Procedure \textproc{E2} is applied to graphs in $B_{m-1, n}$ and treats an input graph as~$G+e_1$ in operation~D3, as expressed in Theorem \ref{maintheorem-1}. It adds all possible edges with a vertex in common to the edge added by \textproc{E1} to yield a graph $G+\{e_1, e_2\}$. This is the second step in operation D3 as expressed in~Theorem~\ref{maintheorem-1}. Cycles in these graphs are also constructed using \textproc{ApplyAddEdge}. The results, after checking certificates, are added to $C_{m, n}$.

\medskip

\item Procedure \textproc{C1} is applied to graphs in $B_{m-1, n-1}$ and treats an input graph as $G + e$ as defined in  operations D1 and D2, as expressed in Theorem \ref{maintheorem-1}. It generates two splits $(G+e) \circ f$ for each input graph, one for each of the vertices incident to the edge added by \textproc{E1}. This is the second step in operations D1 and D2, and it is the final step in D1. It uses \textproc{ApplySubdivideEdge} and \textproc{ApplyFlipEdge} to propagate cycles through the vertex split. The~results, after checking certificates, are added to $A^{(1)}_{m, n}$. This procedure only produces splits for graphs for which the original set of vertices and edges is 3-compatible, and as a result it yields only minimally 3-connected graphs.

\medskip

\item Procedure \textproc{C2} is applied to graphs in $A^{(1)}_{m-1, n-1}$ and treats an input graph as $(G + e) \circ f_c$ as defined in operation D2, as expressed in Theorem \ref{maintheorem-1}. It generates splits $(G + e) \circ \{f_c, f_b\}$ of the remaining un-split vertex incident to the edge added by \textproc{E1}. This is the third step of operation D2 when the new vertex is incident with $e$; otherwise it comprises another application of D1. It~uses  \textproc{ApplySubdivideEdge} and \textproc{ApplyFlipEdge} to propagate cycles through the vertex split.  The~results, after checking certificates, are added to $A^{(2)}_{m, n}$. This procedure only produces splits for 3-compatible input sets, and as a result it yields only minimally 3-connected graphs.

\medskip

\item Procedure \textproc{C3} is applied to graphs in $C_{m-1, n-1}$ and treats an input graph as $G + \{e_1, e_2\}$ as defined in  operation D3 as expressed in Theorem \ref{maintheorem-1}. For each input graph, it generates one vertex split $(G + \{e_1, e_2\}) \circ f$ of the vertex common to the edges added by \textproc{E1} and \textproc{E2}. It~uses \textproc{ApplySubdivideEdge} and \textproc{ApplyFlipEdge} to propagate cycles through the vertex split.  The~results, after checking certificates, are added to $A^{(3)}_{m, n}$. This procedure only produces splits for 3-compatible input sets, and as a result it yields only minimally 3-connected graphs.

\end{itemize}

\medskip

\noindent While \textproc{C1}, \textproc{C2}, and \textproc{C3} produce only minimally $3$-connected graphs, they may produce different graphs that are isomorphic to one another. We use Brendan McKay's nauty to generate a canonical label for each graph produced, so that only pairwise non-isomorphic sets of minimally $3$-connected graphs are ultimately output.

\medskip

The following procedures are defined informally:

\begin{itemize}

\item \textproc{AddEdge}($G, u, v$)-Given a graph $G$ and a pair of vertices $u$ and $v$ in $G$, this procedure returns a~graph $G'$ formed from $G$ by adding an edge connecting $u$ and $v$. When it is used in the procedures in this section, we also use \textproc{ApplyAddEdge} immediately afterwards, which computes the cycles of the graph with the added edge. The complexity of \textproc{AddEdge} is $\mathcal O(n^2)$ because the set of edges of $G$ must be copied to form the set of edges of $G'$.

\medskip

\item \textproc{SplitVertex}($G, v, u, w$) - Given a graph $G$, a vertex $v$ and two edges $v u$ and $v w$, this procedure returns a graph $G'$ formed from $G$ by adding a vertex $v'$, adding an edge connecting $v$ and $v'$, and~replacing the edges $v u$ and $v w$ with edges $v' u$ and $v' w$. When it is used in the procedures in this section, we also use \textproc{ApplySubdivideEdge} and \textproc{ApplyFlipEdge}, which compute the cycles of the graph with the split vertex. The complexity of \textproc{SplitVertex} is $\mathcal O(n^2)$, again because a copy of the graph must be produced.

\medskip

\item \textproc{NoChordingPaths}($\mathcal C, G, P, X$) - Given the set $\mathcal C$ of cycles of a graph $G$, a set $P$ of pairs of vertices and another set $X$ of edges, this procedure determines whether there are any chording paths connecting pairs of vertices in $P$ in $G \backslash X$. Its complexity is $\mathcal O(|\mathcal C| n^3)$, as it requires all simple paths between two vertices to be enumerated, which is $\mathcal O(n^2)$. This function relies on \textproc{HasChordingPath} as defined in Section \ref{correctness-section}.

\end{itemize}

\medskip

The rest of this subsection contains a detailed description and pseudocode for procedures \textproc{E1}, \textproc{E2}, \textproc{C1}, \textproc{C2} and \textproc{C3}. The worst-case complexity for any individual procedure in this process is the complexity of \textproc{C2}: $\mathcal O(|\mathcal C|^2 n^3)$.
\medskip

Procedure \textproc{E1} is responsible for implementing the first step of operations D1, D2, and D3. It~generates all single-edge additions of an input graph $G$, using \textproc{ApplyAddEdge} to propagate the list of cycles. Its complexity is $\mathcal{O}(|\mathcal C|n^3)$, as it requires each pair of vertices of $G$ to be checked, and~for each non-adjacent pair \textproc{ApplyAddEdge} is used to propagate cycles. Pseudocode is shown in Algorithm \ref{alg03}.

\medskip

\begin{algorithm}
\caption{First edge addition procedure}\label{alg03}
{\scriptsize
\begin{algorithmic}[1]
\Procedure{E1}{$G$, $\mathcal C$}\vspace{3pt}
\State $S \gets \phi$\vspace{3pt}
\For{$u \in V(G)$}\vspace{3pt}
	\For{$v \in V(G) \backslash u$}\vspace{3pt}
		\If{$u v \notin E(G)$}\vspace{3pt}
			\State $G' \gets \Call{AddEdge}{G, u, v}$\vspace{3pt}
			\State $\mathcal C' \gets \Call{ApplyAddEdge}{\mathcal C, u, v}$\vspace{3pt}
			\State $S \gets S \cup \{(G', \mathcal C', u v)\}$\vspace{3pt}
		\EndIf\vspace{3pt}
	\EndFor\vspace{3pt}
\EndFor\vspace{3pt}
\State \textbf{return} $S$
\EndProcedure
\end{algorithmic}
}
\end{algorithm}
\vspace{6pt}
\medskip

Procedure \textproc{E2} is responsible for implementing the second step of operation D3. It also generates single-edge additions of an input graph, but under a certain condition. Specifically, given an input graph $G = H + e_1$ with cycles $\mathcal C$, as produced by \textproc{E1}, \textproc{E2} produces all graphs $H + \{e_1, e_2\}$, where the new edge $e_2$ is adjacent to $e_1$. Its complexity is $\mathcal O(|\mathcal C|n^2)$, as \textproc{ApplyAddEdge} is used every time a new graph is generated, and each vertex is checked for eligibility. Pseudocode is shown in Algorithm \ref{alg04}.
\vspace{6pt}
\medskip

\begin{algorithm}
\caption{Second edge addition procedure}\label{alg04}
{\scriptsize
\begin{algorithmic}[1]
\Procedure{E2}{$G$, $u$, $v$, $\mathcal C$}\vspace{3pt}
\State $S \gets \phi$\vspace{3pt}
\For{$w \in V(G)$}\vspace{3pt}
	\If{$wu \notin E(G)$}\vspace{3pt}
		\State $G' \gets \Call{AddEdge}{G, u, w}$\vspace{3pt}
		\State $\mathcal C' \gets \Call{ApplyAddEdge}{\mathcal C, u, w}$\vspace{3pt}
		\State $S \gets S \cup \{(G', \mathcal C', uw)\}$\vspace{3pt}
	\EndIf\vspace{3pt}
	
	\If{$wv \notin E(G)$}\vspace{3pt}
		\State $G' \gets \Call{AddEdge}{G, w, v}$\vspace{3pt}
		\State $\mathcal C' \gets \Call{ApplyAddEdge}{\mathcal C, w, v}$\vspace{3pt}
		\State $S \gets S \cup \{(G', \mathcal C', wv)\}$\vspace{3pt}
	\EndIf\vspace{3pt}
\EndFor\vspace{3pt}
\State \textbf{return} $S$\vspace{3pt}
\EndProcedure\vspace{3pt}
\end{algorithmic}
}
\end{algorithm}

\medskip

Procedure \textproc{C1} is responsible for implementing the second step of operations D1 and D2. These~steps are illustrated in Figures \ref{BG-1-equivalent} and \ref{BG-2-equivalent}, respectively, though a bit of bookkeeping is required to see how \textproc{C1} corresponds to those operations. \textproc{C1} starts with a graph $G = H + e$ generated by \textproc{E1}; let $bc$ denote the added edge. Observe that the chording path checks are made in $H$, which is $G \backslash bc$.

\medskip

First, for any vertex $a \in N(b) \backslash c$ in $G$, where there are no chording $c-a$ or $b-c$ paths in $G \backslash \{bc, ba\} = H \backslash ba$, we split $b$ to add a new vertex $x$ adjacent to $a$, $b$, and $c$. This is the same as the second step illustrated in Figure \ref{BG-1-equivalent} with $c$, $b$, $a$, and $x$ in \textproc{C1} corresponding to $x$, $a$, $b$, and $y$ in the figure, respectively. It is also the same as the second step illustrated in Figure \ref{BG-2-equivalent}, with $c$, $b$, $a$, and $x$ in \textproc{C1} corresponding to $b$, $c$, $d$, and $y$ in the figure, respectively.

\medskip

Second, we must consider splits of the other end vertex of the newly added edge $e$, namely $c$. For~any vertex $d \in N(c) \backslash b$ in $G$ where there are no chording $b-d$ or $b-c$ paths in $G \backslash \{bc, cd\} = H \backslash cd$, we split $c$ to add a new vertex $y$ adjacent to $b$, $c$, and $d$. This is the same as the second step illustrated in Figure \ref{BG-1-equivalent} with $b$, $c$, $d$, and $y$ in \textproc{C1} corresponding to $x$, $a$, $b$, and $y$ in the figure, respectively. It is also the same as the second step illustrated in Figure \ref{BG-2-equivalent}, with $b$, $c$, $d$, and $y$ in \textproc{C1} corresponding to $b$, $c$, $d$, and $y$ in the figure, respectively.

\medskip

Since \textproc{C1} must make $\mathcal O(n)$ calls to \textproc{ApplyFlipEdge}, where $n = |V(G)|$, its complexity is $\mathcal O(|\mathcal C|^2 n^2)$. Pseudocode is shown in Algorithm \ref{alg05}.

\medskip

\begin{algorithm}
\caption{First vertex split procedure}\label{alg05}
{\scriptsize
\begin{algorithmic}[1]
\Procedure{C1}{$G$, $b$, $c$, $\mathcal C$}\vspace{3pt}
\State $S \gets \phi$\vspace{3pt}
\For{$a \in N(b) \backslash c$} \Comment{Split $b$}\vspace{3pt}
	\If{$\Call{NoChordingPaths}{\mathcal C, G, \{(c, a), (b, c)\}, \{b c, b a\}}$}\vspace{3pt}
		\State $G', x \gets \Call{SplitVertex}{G, b, c, a}$\vspace{3pt}
		\State $\mathcal C' \gets \Call{ApplySubdivideEdge}{\mathcal C, b, c, x}$\vspace{3pt}
		\State $\mathcal C' \gets \Call{ApplyFlipEdge}{\mathcal C', a, b, x}$\vspace{3pt}
		\State $S \gets S \cup \{(G', \mathcal C', b, c, a, x)\}$\vspace{3pt}
	\EndIf\vspace{3pt}
\EndFor\vspace{3pt}

\For{$d \in N(c) \backslash b$} \Comment{Split $c$}\vspace{3pt}
	\If{$\Call{NoChordingPaths}{\mathcal C, G, \{(b, d), (c, b)\}, \{b c, c d\}}$}\vspace{3pt}
		\State $G', y \gets \Call{SplitVertex}{G, c, b, d}$\vspace{3pt}
		\State $\mathcal C' \gets \Call{ApplySubdivideEdge}{\mathcal C, c, b, y}$\vspace{3pt}
		\State $\mathcal C' \gets \Call{ApplyFlipEdge}{\mathcal C', d, c, y}$\vspace{3pt}
		\State $S \gets S \cup \{(G', \mathcal C', c, b, d, y)\}$\vspace{3pt}
	\EndIf\vspace{3pt}
\EndFor

\State \textbf{return} $S$\vspace{3pt}
\EndProcedure\vspace{3pt}
\end{algorithmic}
}
\end{algorithm}

\medskip

Procedure \textproc{C2} is responsible for implementing the third step in operation D2, as illustrated in Figure \ref{BG-2-equivalent}. It starts with a graph $G = (H + e) \circ f$ generated by \textproc{C1}; we denote $e'$ and $f$ shown in the figure as $yb$ and $yc$, respectively. $d$ represents the third vertex that becomes adjacent to the new vertex in \textproc{C1}, so $d$ and $y$ are also adjacent.

\medskip

First, for any vertex $a$ adjacent to $b$ other than $c$, $d$, or $y$, for which there are no $c-a$, $c-b$, $d-b$, or $d-a$ chording paths in $G \backslash \{ab, by, cy, dy\}$, we split $b$ to add a new vertex $x$ adjacent to $b$, $a$ and $y$. This is the same as the third step illustrated in Figure \ref{BG-2-equivalent}.

\medskip

Second, for any pair of vertices $a$ and $k$ adjacent to $b$ other than $c$, $d$, or $y$, and for which there are no $k-a$ or $k-b$ chording paths in $G \backslash \{ab, by, cy, dy\}$, we split $b$ to add a new vertex $x$ adjacent to $b$, $a$~and $k$ (leaving $y$ adjacent to $b$, unlike in the first step). 

\medskip

Since \textproc{C2} must make $\mathcal O(n^2)$ calls to \textproc{ApplyFlipEdge}, where $n = |V(G)|$, its complexity is $\mathcal O(|\mathcal C|^2 n^3)$. Pseudocode is shown in Algorithm \ref{alg06}.

\medskip

\begin{algorithm}
\caption{Second vertex split procedure}\label{alg06}
{\scriptsize
\begin{algorithmic}[1]
\Procedure{C2}{$G, \mathcal C, c, b, d, y$}\vspace{3pt}
\State $S \gets \phi$\vspace{3pt}
\For{$a \in N(b) \backslash \{c, d, y\}$} \Comment{Final step of Operation (c)}\vspace{3pt}
	\If{$\Call{NoChordingPaths}{\mathcal C, G, \{ (c,a), (c,b), (d,b), (d,a) \}, \{ab, by, cy, dy\}}$}\vspace{3pt}
		\State $G', x \gets \Call{SplitVertex}{G, b, y, a}$\vspace{3pt}
		\State $\mathcal C' \gets \Call{ApplySubdivideEdge}{\mathcal C, b, a, x}$\vspace{3pt}
		\State $\mathcal C'' \gets \Call{ApplyFlipEdge}{\mathcal C', y, b, x}$\vspace{3pt}
		\State $S \gets S \cup \{(G', \mathcal C'')\}$\vspace{3pt}
	\EndIf\vspace{3pt}
\EndFor\vspace{3pt}
\For{$a, k \in N(b) \backslash \{c, d, y\}$} \Comment{Final step of Operation (d)}\vspace{3pt}
	\If{$\Call{NoChordingPaths}{\mathcal C, G, \{ (k,a), (k,b) \}, \{ab, by, cy, dy\}}$}\vspace{3pt}
		\State $G', x \gets \Call{SplitVertex}{G, b, k, a}$\vspace{3pt}
		\State $\mathcal C' \gets \Call{ApplySubdivideEdge}{\mathcal C, b, a, x}$\vspace{3pt}
		\State $\mathcal C'' \gets \Call{ApplyFlipEdge}{\mathcal C', k, b, x}$\vspace{3pt}
		\State $S \gets S \cup \{(G', \mathcal C'')\}$\vspace{3pt}
	\EndIf\vspace{3pt}
\EndFor\vspace{3pt}
\State \textbf{return} $S$\vspace{3pt}
\EndProcedure\vspace{3pt}
\end{algorithmic}
}
\end{algorithm}

\medskip

Procedure \textproc{C3} is responsible for implementing the third step in operation D3, as illustrated in Figure \ref{D-3-equivalent}. It starts with a graph $G = H + \{e_1, e_2\}$ generated by \textproc{E2}, where $e_1=xy$ and $e_2=xz$ are two incident edges. A single new graph is generated in which $x$ is split to add a new vertex $w$ adjacent to $x$, $y$ and $z$, if there are no $x-y$, $x-z$, or $y-z$ chording paths in $G \backslash \{xy, xz\} = H$. Because \textproc{C3} makes one call to \textproc{ApplyFlipEdge}, its complexity is $\mathcal O(|\mathcal C|^2 n)$. Pseudocode is shown in Algorithm \ref{alg07}.

\medskip

\begin{algorithm}
\caption{Third vertex split procedure}\label{alg07}
{\scriptsize
\begin{algorithmic}[1]
\Procedure{C3}{$G, \mathcal C, x, y, z$}\vspace{3pt}
\State $S \gets \phi$\vspace{3pt}

\If{$\Call{NoChordingPaths}{\mathcal C, G, \{(x, y), (x, z), (y, z)\}, \{xy, xz\}}$}\vspace{3pt}
	\State $G', w \gets \Call{SplitVertex}{G, x, y, z}$\vspace{3pt}
	\State $\mathcal C' \gets \Call{ApplySubdivideEdge}{\mathcal C, x, z, w}$\vspace{3pt}
	\State $\mathcal C' \gets \Call{ApplyFlipEdge}{\mathcal C', y, x, w}$\vspace{3pt}
	\State $S \gets S \cup \{(G', \mathcal C')\}$	\vspace{3pt}
\EndIf\vspace{3pt}

\State \textbf{return} $S$\vspace{3pt}
\EndProcedure\vspace{3pt}
\end{algorithmic}
}
\end{algorithm}

%%%%%%%%%%%%%%%%%%%%%%%%%%%%%%%

\section{\bf The Algorithm Is Exhaustive}\label{exhaustive-section}

\bigskip

Theorem \ref{DawesTheorem} and Theorem \ref{DawesTheorem2} (Dawes' results) state that, if $G$ is a minimally 3-connected graph and $G'$ is obtained from $G$ by applying one of the operations D1, D2, and D3 to a set $S$ of vertices and edges, then $G'$ is minimally 3-connected if and only if $S$ is 3-compatible, and also that any minimally 3-connected graph other than $W_3$ can be obtained from a smaller minimally 3-connected graph by applying D1, D2, or D3 to a 3-compatible set.  This shows that application of these operations to 3-compatible sets of edges and vertices in minimally 3-connected graphs, starting~with $W_3$, will~exhaustively generate all such graphs. However, as indicated in Theorem \ref{cycle-prop-theorem}, in order to maintain the list of cycles of each generated graph, we must express these operations in terms of edge additions and vertex splits. 

\medskip

\begin{figure}[H]
\centering
\includegraphics[width=3.4in]{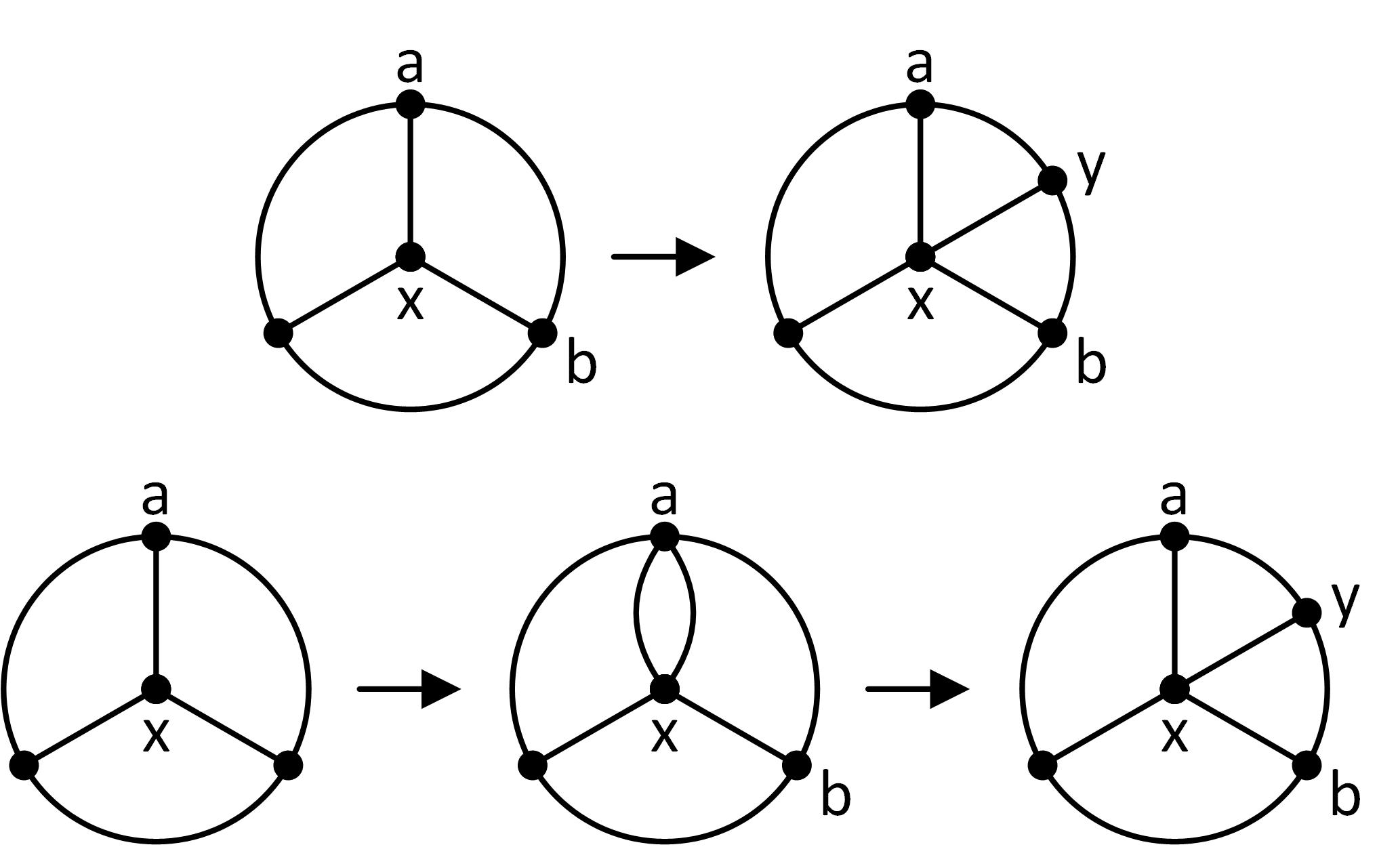}
\caption{Operation D1 applied to $W_3$. \label{d1-applied-to-k4}}
\end{figure}

Consider the graph $W_3$ itself, as shown in Figure \ref{d1-applied-to-k4}. Observe that $\{x, ab\}$ is a 3-compatible set because there are clearly no chording $xa$- or $xb$-paths in $W_3 \backslash ab$, so we may apply D1 to produce another minimally 3-connected graph, which is actually $W_4$ as shown in the figure. However, since~there are already edges $xa$ and $xb$ in the graph, if we are to apply our step-by-step procedure to accomplish the same thing, we will be required to add a parallel edge. We would like to avoid this, and we can accomplish that by beginning with the prism graph instead of $W_3$. We are now ready to prove the third main result in this paper.

\begin{theorem}\label{exhaustive-theorem-2} Let $G$ be a simple minimally $3$-connected graph. Then one of the following statements is true:
\begin{enumerate}
\item $G \cong W_{n-1}$ for $n \ge 5$ and $G$ can be obtained from $W_{n-2}$ by applying operation D1 to the spoke vertex $x$ and a rim edge $ab$;
\item $G \cong K_{3, n-3}$ for $n \ge 7$ and $G$ can be obtained from $K_{3, n-4}$ by applying operation D3 to the $3$ vertices in the smaller class; or
\item $G$ has a prism minor, for $n \ge 7$, and $G$ can be obtained from a smaller minimally $3$-connected graph $G'$ with a prism minor, where $|E(G)| - |E(G')| \le 3$, using operation D1, D2, or D3.
\end{enumerate}
\end{theorem}

\begin{proof} Theorem \ref{DiracTheorem} characterizes the $3$-connected graphs without a prism minor. Of these, the only minimally $3$-connected ones are $W_{n-1}$ for $n \ge 4$ and $K_{3, n-3}$ for $n \ge 6$. Observe that for $n\ge 5$, $W_{n-1}\backslash e/f = W_{n-2}$, where $e$ is a spoke and $f$ is a rim edge, such that $e, f$ are incident to a degree 3~vertex. Therefore $W_{n-1}$ can be obtained from $W_{n-2}$ by applying operation D1 to the spoke vertex $x$ and a rim edge $ab$. The set $\{x, ab\}$ is 3-compatible because any chording edge of a cycle in $W_{n-2} \backslash ab$ would have to be a spoke edge, and since all rim edges have degree three the chording edge cannot be extended into a $xa$- or $xb$-path. 

\medskip

Observe that, for $n\ge 7$, $K_{3, n-3}-w = K_{3, n-4}$, where $w$ is a degree 3 vertex. Therefore, $K_{3, n-3}$ can be obtained from a smaller minimally 3-connected graph of the same family by applying operation D3 to the three vertices in the smaller class. The set of three vertices is 3-compatible because the degree of each vertex in the larger class is exactly 3, so that any chording edge cannot be extended into a~chording path connecting vertices in the smaller class, as illustrated in Figure \ref{k-3-n-3}.

\begin{figure}[H]
\centering
\includegraphics[width=0.6in]{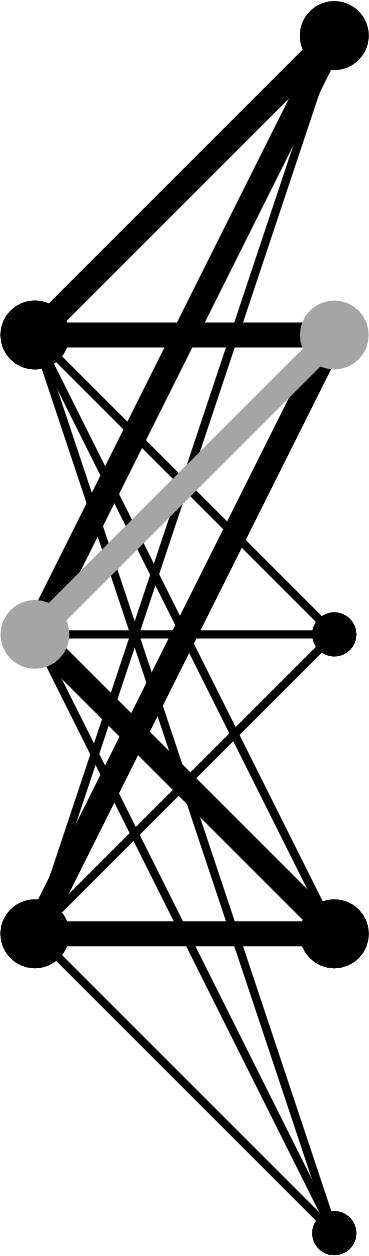}
\caption{$K_{3, 5}$ with a cycle highlighted in black and a chording edge highlighted in gray.\label{k-3-n-3}}
\end{figure}

If $G$ has a prism minor, by Theorem \ref{strong-splitter-theorem}, with the prism graph as $H$, $G$ can be obtained from a~3-connected graph with $n-1$ vertices and $m-1$ edges via an edge addition and a vertex split, from~a~graph with $n-2$ vertices and $m-3$ edges via two edge additions and a vertex split, or from a graph with $n-1$ vertices and $m-3$ edges via an edge addition and two vertex splits; that is, by operation D1, D2, or D3, respectively, as expressed in Theorem \ref{maintheorem-1}.
By Theorem \ref{DawesTheorem2}, all minimally 3-connected graphs can be obtained from smaller minimally 3-connected graphs by applying these operations to 3-compatible sets. 
\end{proof}
 \medskip

We constructed all non-isomorphic minimally $3$-connected graphs up to 12 vertices using a Python implementation of these procedures. The total number of minimally $3$-connected graphs for 4 through 12 vertices is published in the Online Encyclopedia of Integer Sequences. Table \ref{graphs-table} below lists these values. These numbers helped confirm the accuracy of our method and procedures.

\begin{table}[H]

\caption{Minimally 3-connected graphs from \texttt{oeis.org}. \label{graphs-table}}
\centering
\begin{tabular}{ccccc} \hline
\bf   \boldmath{$n$}  & \bf Minimally 3-Connected Graphs (A199676) \\  \hline
6   &    3   \\ \hline
7     &     5\\ \hline 
8   &     18\\ \hline
9    &     57\\ \hline
10   &   285\\ \hline
11   &   1513\\ \hline
12 & 9824\\ \hline
\end{tabular}
\end{table}

The number of non-isomorphic 3-connected cubic graphs of size $n$, where $n$ is even, is published in the Online Encyclopedia of Integer Sequences as sequence A204198. This sequence only goes up to $n=14$. We were able to obtain the set of 3-connected cubic graphs up to 20 vertices as shown in~Table~\ref{cubic-table}.

\begin{table}[H]

\caption{3-connected cubic graphs. \label{cubic-table}}
\centering
\begin{tabular}{ccccc} \hline
\bf  \boldmath{$n$} & \bf{3-Connected Cubic Graphs} \\  \hline
8 & 4   \\ \hline
10 & 14 \\ \hline 
12 & 57 \\ \hline 
14 & 341 \\ \hline
16 & 2828 \\ \hline 
18 & 30,468 \\ \hline
20 & 396,150 \\ \hline
\end{tabular}
\end{table}
\medskip

All of the minimally 3-connected graphs generated were validated using a separate routine based on the Python iGraph ({\texttt{https://igraph.org/python/}}) {\tt vertex\_disjoint\_paths} method, in order to verify that each graph was $3$-connected and that all single edge-deletions of the graph were not. The~overall number of generated graphs was checked against the published sequence on OEIS.  

\medskip

The $3$-connected cubic graphs were verified to be $3$-connected using a similar procedure, and~overall numbers for up to 14 vertices were checked against the published sequence on OEIS.

\medskip

The minimally $3$-connected graphs were generated in 31 h on a PC with an Intel Core I5-4460 CPU at 3.2 GHz and 16 Gb of RAM. The 3-connected cubic graphs were generated on the same machine in five hours.

\medskip

The algorithm’s running speed could probably be reduced by running parallel instances, either on a larger machine or in a distributed computing environment. MapReduce, or a similar programming model, would need to be used to aggregate generated graph certificates and remove duplicates. It is also possible that a technique similar to the canonical construction paths described by Brinkmann, Goedgebeur and McKay  [\ref{Brinkmannetal2011}] could be used to reduce the number of redundant graphs generated.

\medskip

Even with the implementation of techniques to propagate cycles, the slowest part of the algorithm is the procedure that checks for chording paths. It may be possible to improve the worst-case performance of the cycle propagation and chording path checking algorithms through appropriate indexing of cycles.

\medskip

The code, instructions, and output files for our implementation are all available at \texttt{https://github.com/rkingan/m3c}. The output files have been converted from the format used by the program, which~also stores each graph's history and list of cycles, to the standard graph6 format, so that they can be used by other researchers.

\bigskip
\noindent \textbf{Acknowledgements:} The authors would like to thank the referees and editor for their valuable comments which helped to improve the manuscript.

 \bigskip

\section*{\bf References}

\bigskip

\begin{enumerate}
\item \label{BarnetteGrunbaum1969} D. W. Barnette and B. Gr\"unbaum (1969). On Steinitz's theorem concerning convex 3-polytopes and on some properties of planar graphs, {\it The Many Facets of Graph Theory. Lecture Notes in Mathematics}, vol 110, G. Chartrand  and  S. F. Kapoor (eds),  Springer, Berlin, Heidelberg, 27--40.

\item \label{Brinkmannetal2011} G. Brinkmann, J. Goedgebeur and B.D. Mckay (2011). Generation of Cubic graphs. \textit{Discrete Mathematics and Theoretical Computer Science, DMTCS}, \textbf{13(2)}, 69--79.

\item \label{cormen} T. H. Cormen, C. E. Leiserson, R. L. Rivest, and C. Stein (1989). Introduction to Algorithms, Third Edition, 	MIT Press. 
 
\item  \label{Dawes1986}  R. W. Dawes (1986). Minimally 3-connected graphs, {\it J. Combin. Theory Ser. B} {\bf 40}, 159-168. 
 
\item  \label{Dirac1963}  G. A. Dirac (1963). Some results concerning the structure of graphs, {\it Canad. Math. Bull.} {\bf 6}, 183-210. 

\item  \label{Halin1969} R. Halin (1969). Untersuchungen uber minimale n-fach zusammenhangende graphen, {\it Math. Ann} 182 (1969), 175--188.
 
\item \label{Hopcroft1973} J. E. Hopcroft and R. E. Tarjan (1973). Dividing a graph into triconnected components.
{\it SIAM J. Comput.}, {\bf 2(3)}, 135 – 158.

\item  \label{KinganLemos2014}	S. R. Kingan  and M. Lemos  (2014). Strong Splitter Theorem, {\it Annals of Combinatorics} {\bf 18-1}, 111-116.

\item \label{McKay1981} McKay, B.D. Practical Graph Isomorphism. {\it Congr. Numer.} {\bf 1981}, {\em 30}, 45--87.

\item  \label{Oxley2011} J. G. Oxley (2011).  {\it Matroid Theory}, Second Edition, Oxford University Press (2011),  New York.
 
\item \label{Schmidt2011} J. M. Schmidt  (2011). {\it Structure and constructions of $3$-connected graphs.} (Dissertation) Free University of Berlin.

\item \label{Tutte1961}  W. T. Tutte (1961). A theory of 3-connected graphs, {\it Indag. Math} {\bf 23}, 441-455.

\item \label{Tutte1967} W.T. Tutte (1967). {\it Connectivity in Graphs}. Toronto University Press, Toronto.

\end{enumerate}

%%%%%%%%%%%%%%%%%%%%%%%%%%%%%%%%%%%%%%%%%%%%%%%%%%%

\end{document}